\documentclass[12pt]{article}
\usepackage{amscd, amssymb, amsmath}
\begin{document}
\renewcommand{\thesubsection}{\arabic{subsection}}
\newenvironment{proof}{{\bf Proof}:}{\vskip 5mm }
\newenvironment{rem}{{\bf Remark}:}{\vskip 5mm }
\newenvironment{remarks}{{\bf Remarks}:\begin{enumerate}}{\end{enumerate}}
\newenvironment{examples}{{\bf Examples}:\begin{enumerate}}{\end{enumerate}}  
\newtheorem{proposition}{Proposition}[subsection]
\newtheorem{lemma}[proposition]{Lemma}
\newtheorem{definition}[proposition]{Definition}
\newtheorem{theorem}[proposition]{Theorem}
\newtheorem{cor}[proposition]{Corollary}
\newtheorem{conjecture}{Conjecture}
\newtheorem{pretheorem}[proposition]{Pretheorem}
\newtheorem{hypothesis}[proposition]{Hypothesis}
\newtheorem{example}[proposition]{Example}
\newtheorem{remark}[proposition]{Remark}
\newtheorem{ex}[proposition]{Exercise}
\newtheorem{cond}[proposition]{Conditions}
\newtheorem{constr}[proposition]{Construction}
\newcommand{\llabel}[1]{\label{#1}}
\newcommand{\comment}[1]{}
\newcommand{\sr}{\rightarrow}
\newcommand{\dw}{\downarrow}
\newcommand{\bdl}{\bar{\Delta}}
\newcommand{\cc}{{\bf C}}
\newcommand{\zz}{{\bf Z\rm}}
\newcommand{\zq}{{\bf Z}_{qfh}}
\newcommand{\nn}{{\bf N\rm}}
\newcommand{\qq}{{\bf Q\rm}}
\newcommand{\nq}{{\bf N}_{qfh}}
\newcommand{\oo}{\otimes}
\newcommand{\uu}{\underline}
\newcommand{\ih}{\uu{Hom}}
\newcommand{\af}{{\bf A}^1}
\newcommand{\dsr}{\stackrel{\sr}{\scriptstyle\sr}}
\newcommand{\PP}{$P_{\infty}$}
\newcommand{\DD}{D}
\newcommand{\tp}{\tilde{D}}
\newcommand{\HH}{$H_{\infty}$}
\newcommand{\ii}{\stackrel{\scriptstyle\sim}{\sr}}
\newcommand{\BB}{_{\bullet}}

\begin{center}
{\Large\bf Reduced power operations in motivic
cohomology}\footnote{Supported by the NSF grants 
DMS-97-29992 and DMS-9901219 and The Ambrose Monell
Foundation}$^,$\footnote{School of Mathematics, Institute for Advanced
Study, Princeton NJ, USA. e-mail: vladimir@ias.edu}\\
\vskip 4mm
{\large\bf by Vladimir Voevodsky}\\
\vskip 3mm
{\large May 2001}
\end{center}
\vskip 4mm
\tableofcontents

\numberwithin{equation}{subsection}

\subsection{Introduction}
In this paper we consider cohomological operations in the motivic
cohomology of smooth simplicial schemes over a field $k$. For the most
part we work with cohomology with coefficients in $\zz/l$ where $l$ is
a prime different from the characteristic of $k$.  We construct the
reduced power operations 
$$P^i:H^{*,*}(-,\zz/l)\sr H^{*+2i(l-1),*+i(l-1)}(-,\zz/l)$$
and prove the motivic analogs of the Cartan formula and the Adem
relations. We also describe the subalgebra in the algebra of all
(bistable) operations in the motivic cohomology with
$\zz/l$-coefficients generated by operations $P^i$, the Bockstein
homomorphism and the multiplication by the motivic cohomology classes
of $Spec(k)$. For odd $l$ this algebra is isomorphic to the twisted
tensor product of the usual topological Steenrod algebra with the
motivic cohomology ring of the point $H^{*,*}$ with respect to the
action of the motivic Steenrod operations on $H^{*,*}$. For $l=2$, the
situation is more complicated since the motivic Adem relations involve
nontrivial coefficients from $H^{*,*}$.

To construct the reduced power operations we follow the approach of
\cite{SE} where one first defines the total power operation and then
uses the computation of the cohomology of the product of a space with
the classifying space of the symmetric group $S_l$ to obtain the
individual power operations. We also use the ideas of \cite{SE} for
the proofs of the relations between power operations and the
Bockstein homomorphism, the Cartan formula and the Adem
relations. Our construction of the total power operations is not
directly similar to any of the standard topological
constructions\footnote{The construction of the total power operations
given here is also slightly different from the construction given in
my unpublished notes on the operations written in 1996.}. One would get
its direct topological analog if one unfolded, to the space-level, the
description of power operations in terms of $E_{\infty}$-structure on
the Eilenberg-MacLane spectra. In particular, the Thom isomorphism in
the motivic cohomology and the Euler classes of vector bundles figure
prominently in most of our computations.

Several important results on operations in motivic cohomology are not
discussed in this paper and will be proved in a sequel. We do not show
that the operations $P^i$ are unique. We do not show that the operations $P^i$
and the Bockstein homomorphism generate all bistable cohomological
operations. Finally, we do not prove that bistable cohomological
operations coincide with the endomorphisms of the Eilenberg-Maclane
spectrum in the stable category. 

Modulo this identification of bistable operations with the
endomorphisms of $H_{\zz/l}$, the present paper contains proofs of the
following results used in \cite{MC} for the proof of the Milnor
conjecture: \cite[Theorem 3.14 p.31]{MC} is Propositions \ref{neg} and
\ref{zer}; \cite[Theorem 3.16 p.32]{MC} is Theorem \ref{zero}, Lemma
\ref{power}, Proposition \ref{cart} and Lemma \ref{l2.6.8}. As Example
\ref{incor} shows, the inductive construction used in \cite{MC} to define
the operations $Q_i$ is incorrect unless $k$ contains a square root of
$-1$. Instead we define them in a different way and prove \cite[Theorem
3.17 p.32]{MC} in Proposition \ref{bstrd} and Proposition
\ref{com1}. This paper also contains all the results about
cohomological operations necessary for the proof of the Milnor
Conjecture given in \cite{MCnew}.

The first draft of this paper was written in April 1996 i.e., exactly
five years ago. During these years I discussed problems related to
operations in motivic cohomology with a lot of people and I am
greatful to all of them for these conversations. I would like to
especially thank Mike Hopkins, Fabien Morel, Charles Rezk and
Alexander Vishik. 

This paper was written while I was a member of the Institute for
Advanced Study in Princeton and, part of the time, an employee of the
Clay Mathematics Institute. I am very grateful to both institutions
for their support.

\subsection{Motivic cohomology and operations}
For any $p,q\in\zz$ and any abelian group $A$ the motivic cohomology
of a smooth scheme $X$ are defined as hypercohomology
$$H^{p,q}(X,A):={\bf H}^p(X_{Nis},A(q))$$
where $A(q)=\zz(q)\oo A$ and $\zz(q)$ is a certain complex of sheaves
of abelian groups on $(Sm/k)_{Nis}$. Let $K(p,q,A)$ be the simplicial
abelian sheaf corresponding to the complex $A(q)[p]$. Considered as a
pointed simplicial sheaf of sets it defines an object of the pointed
motivic homotopy category $H_{\bullet}(k)$. The simplicial sheaves
$K(p,q,A)$ are $\af$-local and for a smooth scheme $X$ one has
$$Hom_{H_{\bullet}(k)}(X_+,K(p,q,A))=H^{p,q}(X,A)$$
i.e. the objects $K(p,q,A)$ represent motivic cohomology on
$H_{\bullet}$. For any pointed simplicial sheaf $F_{\bullet}$ we
define its reduced motivic cohomology by 
$$\tilde{H}^{p,q}(F_{\bullet},A):=Hom_{H_{\bullet}(k)}(F_{\bullet},K(p,q,A))$$
Let $R$ be a commutative ring with unit. Then $R(q)$ are complexes of
sheaves of free $R$-modules and one has multiplication maps
$$R(q)\oo R(q')\sr R(q+q')$$
which define morphisms of pointed sheaves
\begin{equation}
\llabel{mult}
K(p,q,R)\wedge K(p',q',R)\sr K(p+p',q+q',R)
\end{equation}
Recall that for a smooth scheme $X$ over $k$ we denote by
$\zz_{tr}(X)$ the presheaf on $Sm/k$ which takes $U$ to the group of
cycles on $U\times X$ which are finite and equidimensional over
$U$. For an abelian group $A$ we denote by $A_{tr}(X)$ the presheaf
$\zz_{tr}(X)\oo A$.

Let $K_{n,A}$ be the pointed sheaf of sets associated with the
presheaf 
$$K^{pre}_{n,A}:U\mapsto A_{tr}({\bf A}^n)(U)/A_{tr}({\bf
A}^n-\{0\})(U)$$
where the quotient is the quotient in the category of abelian
groups. A section of $K_n^{pre}$ on $U$ is an equivalence class of
cycles on $U\times {\bf A}^n$ with coefficients in $A$ which are
equidimensional and finite over $U$. If $R$ is a commutative ring then
we have the multiplication maps
\begin{equation}
\llabel{exm}
K_{n,R}\wedge K_{m,R}\sr K_{n+m,R}
\end{equation}
which send the pair of cycles $Z$, $Z'$ to the external product cycle
$Z\oo Z'$.
The following result is proved in
\cite{delnotes}.
\begin{theorem}
\llabel{del1}
There are $\af$-weak equivalences $K_{n,R}\sr K(2n,n,R)$ which are
compatible with the multiplication maps.
\end{theorem}
For pointed sheaves $F_{\BB}$, $G_{\BB}$ the morphisms (\ref{mult})
define multiplication maps
\begin{equation}
\llabel{mul1}
H^{*,*}(F_{\BB},R)\oo H^{*,*}(G_{\BB},R)\sr H^{*,*}(F_{\BB}\wedge
G_{\BB},R)
\end{equation}
which we denote, on elements, by $w\oo w'\mapsto w\wedge w'$. If
$G_{\BB}=F_{\BB}$ the composition of  (\ref{mul1}) with the map
defined by the diagonal $F_{\BB}\sr F_{\BB}\wedge F_{\BB}$ defines
multiplication maps
\begin{equation}
\llabel{mul2}
H^{*,*}(F_{\BB},R)\oo H^{*,*}(F_{\BB},R)\sr H^{*,*}(F_{\BB},R)
\end{equation}
which we denote on elements by $w\oo w'\mapsto w w'$.
\begin{theorem}
\llabel{comass}\llabel{l38}
The morphisms (\ref{mul2}) define, for any $F_{\bullet}$, a structure
of an associative $R$-algebra on $\tilde{H}^{*,*}(F_{\bullet},R)$
which is graded commutative with respect to the first grading.
\end{theorem}
\begin{proof}
Standard arguments from homological algebra together with the fact
that complexes $R(q)$ do not have cohomology in dimensions $>2q$ imply
that it is sufficient to check that the multiplication maps 
$$K(2q,q,R)\wedge K(2q',q',R)\sr K(2(q+q'),q+q',R)$$
are associative and commutative in $H_{\bullet}$. By Theorem
\ref{del1} it is sufficient to check that the multiplication maps
(\ref{exm}) are commutative and associative in $H_{\bullet}$. The
associativity condition clearly holds on the level of sheaves. To
prove commutativity we should show that permutation of coordinates on
${\bf A}^n$ acts trivially on $K_{n,R}$. The action of the permutation
group on $K_{n,R}$ extends to an action of $GL_n$. A transposition is
$\af$-homotopic to the automorphism given by $(x_1,\dots,x_n)\mapsto
(-x_1,\dots, x_n)$. It is therefore sufficient to check that this
automorphism is the identity in $H_{\bullet}$.  Consider for
simplicity of notations the case of one variable i.e. the automorphism
$\phi$ of $\zz_{tr}(\af)/\zz_{tr}(\af-\{0\})$ defined by $x\mapsto
-x$. The sheaf $\zz_{tr}(\af)/\zz_{tr}(\af-\{0\})$ is isomorphic to
the sheaf $\zz_{tr}({\bf P}^1)/\zz_{tr}({\bf P}^1-\{0\})$ which is
weakly equivalent to $\zz_{tr}({\bf P}^1)/\zz$ where the embedding
$\zz\sr {\bf P}^1$ corresponds to the point $\infty$. Under this weak
equivalence our automorphism becomes the automorphism of
$\zz_{tr}({\bf P}^1)/\zz$ defined by $z\mapsto -z$. Denote this
automorphism by $\psi$. One can easily see now that to prove that
$\phi$ is identity in $H_{\bullet}$ it is sufficient to construct a
section $h$ of $\zz_{tr}({\bf P}^1)$ on ${\bf P}^1\times\af$ such that
\begin{equation}
\llabel{homot}
h_{{\bf
P}^1\times\{1\}}-h_{{\bf
P}^1\times\{-1\}}=Id-\psi
\end{equation}
Let $((z_0:z_1),t)$ be the coordinates on ${\bf P}^1\times\af$ and
$(x_0:x_1)$ the coordinates on ${\bf P}^1$. Then the cycle of the
closed subscheme in ${\bf P}^1\times\af\times{\bf P}^1$ given by the
equation $z_0x_1^2+tz_1x_0x_1+(t^2-1)z_1x_0^2$ defines a section of
$\zz_{tr}({\bf P}^1)$ on ${\bf P}^1\times\af$ which satisfies
(\ref{homot}).
\end{proof}
Denote by $H^{*,*}$ the ring
$$H^{*,*}(Spec(k),R)=\tilde{H}^{*,*}(S^0,R)$$
Then for any $F_{\BB}$ the multiplication maps (\ref{mul1}) define a
structure of $H^{*,*}$-module on $\tilde{H}^{*,*}(F_{\BB},R)$. Theorem
\ref{comass} immediately implies the following fact.
\begin{cor}
\llabel{ismod}
The multiplication map (\ref{mul2}) factors through an
$H^{*,*}$-module map
\begin{equation}
\llabel{mul3} H^{*,*}(F_{\BB},R)\oo_{H^{*,*}} H^{*,*}(F_{\BB},R)\sr
H^{*,*}(F_{\BB},R)
\end{equation}
\end{cor}
Let $S^1_s$ and $S^1_t$ be the circles in $H_{\bullet}$. We have
canonical classes
$$\sigma_s\in \tilde{H}^{1,0}(S^1_s,R)$$
$$\sigma_t\in \tilde{H}^{1,1}(S^1_t,R)$$
Multiplication with these classes gives us suspension morphisms
\begin{equation}
\llabel{ss}
\tilde{H}^{p,q}(F_{\bullet},R)\sr \tilde{H}^{p+1,q}(F_{\bullet}\wedge S^1_s, R)
\end{equation}
\begin{equation}
\llabel{st}
\tilde{H}^{p,q}(F_{\bullet},R)\sr \tilde{H}^{p+1,q+1}(F_{\bullet}\wedge S^1_t, R)
\end{equation}
\begin{theorem}
\llabel{susp}
The suspension morphisms are isomorphisms.
\end{theorem} 
\begin{proof}
Let $\tilde{\zz}$ be the functor from sheaves of pointed sets to
sheaves of abelian groups which sends a sheaf of sets to the freely
generated sheaf of abelian groups with the distinguished point set to
be zero. Let further $N$ be the normalized chain complex functor from
simplicial abelian sheaves to the complexes of abelian sheaves. Then
for any $F_{\bullet}$ one has
$$\tilde{H}^{p,q}(F_{\BB},A)=Hom_{D}(N\tilde{\zz}(F_{\BB}),A(q)[p])$$
where $D$ is the derived category of complexes of abelian sheaves in
the Nisnevich topology. The fact that (\ref{ss}) is an isomorphism
follows from the fact that $N\zz$ takes smash product to tensor
product (modulo a quasi-isomorphism) and that $N\zz(S^1_s)$ is
quasi-isomorphic to $\zz[1]$. 

Consider the suspension morphism
$$(-)\wedge \sigma_T:\tilde{H}^{p,q}(F_{\BB},A)\sr \tilde{H}^{p+2,q+1}(F_{\BB}\wedge
T,A)$$
given by multiplication with the class $\sigma_T\in \tilde{H}^{2,1}(T)$. Since
(\ref{ss}) is an isomorphism and we have an $\af$-weak equivalence
$S^1_s\wedge S^1_t\sr T=h_{\af}/h_{(\af-\{0\})}$, to show that
(\ref{st}) is an isomorphism it is sufficient to show that $(-)\wedge
\sigma_T$ is an isomorphism. A standard argument allows one to reduce
the problem to the case when $F_{\BB}=(h_U)_+$ for a smooth scheme $U$
over $S$. Open excision implies that $T={\bf P}^1/{\bf A}^1$ and we
get a split short exact sequence
$$0\sr \tilde{H}^{*,*}((h_U)_+\wedge T,A)\sr H{*,*}(U\times
{\bf P}^1, A)\sr 
H^{*,*}(U\times \af, A)\sr 0$$
Consider the morphism of sequences:
\begin{equation}
\llabel{mor}
\begin{CD}
H^{*-2,*-1}(U) @>(-)\wedge \sigma>> H^{*,*}(U\times {\bf
P}^1) @>>> H^{*,*}(U\times\af) \\ @VVV @VVV @VVV
\\ \tilde{H}^{*,*}(U_+\wedge T) @>>>
H^{*,*}(U\times {\bf P}^1) @>>> H^{*,*}(U\times\af) 
\end{CD}
\end{equation}
where $\sigma$ is the restriction of $\sigma_T$ to ${\bf P}^1$, the
first vertical arrow is $(-)\wedge\sigma_T$ and the rest of vertical
arrows are identities. The fact that $(-)\wedge \sigma_T$ is an
isomorphism follows now from Lemma \ref{blo} below.
\end{proof}
\begin{lemma}
\llabel{blo}
The upper sequence in (\ref{mor}) is a short exact sequence.
\end{lemma}
\begin{proof}
By \cite[Cor. 2]{comparison} we have natural isomorphisms
$$H^{p,q}(U,\zz)\sr CH^q(U,2q-p)$$
where the target are Bloch's higher
Chow groups and the proof immediately shows that we have similar
isomorphisms for all groups of coefficients
\begin{equation}
\llabel{mor2}
H^{p,q}(U,A)\sr CH^q(U,2q-p,A)
\end{equation}
Consider the diagram:
\begin{equation}
\label{mor1}
\begin{CD}
H^{p-2,q-1}(U) @>(-)\wedge \sigma>> H^{p,q}(U\times {\bf
P}^1) @>>> H^{p,q}(U\times\af) \\ @VVV @VVV @VVV
\\ CH^{q-1}(U, {\scriptstyle{2q-p}}) @>>>
CH^q(U\times {\bf P}^1,{\scriptstyle{2q-p}}) @>>> CH^{q}(U\times\af,{\scriptstyle{2q-p}}) 
\end{CD}
\end{equation}
where in the lower line the first morphism is given by covariant
functoriality for the closed embedding $U\sr U\times {\bf P}^1$ at the
infinity and the second morphism is given by the contravariant
functoriality for the open embedding $U\times\af\sr U\times{\bf
P}^1$. One verifies using the explicit form of the isomorphism
(\ref{mor2}) that both squares in (\ref{mor1}) commute. We conclude
that the upper line is a short exact sequence since the lower one is a
short exact sequence by Bloch's Localization Theorem \cite{SBloch1}.
\end{proof}

We define a bistable cohomological operation of bidegree $(i,j)$ in
motivic cohomology with coefficients in $R$ as a collection of natural
transformations of functors on $H_{\bullet}$
$$\phi_{p,q}:\tilde{H}^{p,q}(-,R)\sr \tilde{H}^{p+i,q+j}(-,R)$$
which commutes with the suspension morphisms i.e. such that for 
$x\in \tilde{H}^{p,q}$ one has
\begin{equation}
\llabel{bistab}
\begin{array}{l}
\phi_{p+1,q}(x\wedge\sigma_s)=\phi_{p,q}(x)\wedge \sigma_s\\
$$\phi_{p+1,q+1}(x\wedge\sigma_t)=\phi_{p,q}(x)\wedge \sigma_t
\end{array}
\end{equation}
Denote by $\sigma_T$ the canonical element in $\tilde{H}^{2,1}(T,R)$ where
$T=\af/\af-\{0\}$. 
\begin{proposition}
\llabel{redtoone}
There is a bijection between the set of bistable cohomological
operation of bidegree $(i,j)$ and  the collections of natural transformations 
$$\phi_n:\tilde{H}^{2n,n}(-,R)\sr \tilde{H}^{2n+i,n+j}(-,R)$$
given for all $n\ge 0$ such that 
\begin{equation}
\llabel{ob1}
\phi_{n+1}(x\wedge \sigma_T)=\phi_n(X)\wedge \sigma_T
\end{equation}
\end{proposition} 
\begin{proof}
Since $T=S^1_s\wedge S^1_t$ and $\sigma_T=\sigma_s\wedge \sigma_t$, the
restriction of a bistable operation to groups of degree $(2n,n)$
satisfies (\ref{ob1}). On the other hand, for a family $\phi_n$ we
can construct $\phi_{p,q}$ as follows. For $F_{\bullet}$ we have
$$\tilde{H}^{p,q}(F_{\bullet})=\tilde{H}^{p+a+b,q+b}(S^a_t\wedge S^b_s\wedge
F_{\bullet})$$
taking $a=2q-p+b$ and taking $b$ to be greater or equal to $max\{0,-q,
p-2q\}$ we get:
$$\tilde{H}^{p,q}(F_{\bullet})=\tilde{H}^{2(q+b),q+b}(S^a_t\wedge S^b_s\wedge
F_{\bullet})$$
where $a,b,q+b\ge 0$. Using these isomorphisms, the operation $\phi_{q+b}$
defines a map
$$\phi_{p,q}:\tilde{H}^{p,q}(F_{\bullet})\sr \tilde{H}^{p+i,q+j}(F_{\bullet})$$
The condition (\ref{ob1}) implies that this map does not depend on the
choice of $b$ and that the maps $\phi_{p,q}$ all $(p,q)$ satisfy
(\ref{bistab}). 
\end{proof}
Combining Theorem \ref{del1} and Proposition \ref{redtoone} we get.
\begin{proposition}
\llabel{descr}
There is a bijection between the set of bistable cohomological
operations of bidegree $(i,j)$ and collections of motivic cohomology
classes $\phi_n\in \tilde{H}^{2n+i,n+j}(K_{n,R},R)$ such that the restriction
of $\phi_{n+1}$ to $K_{n,R}\wedge T$ is $\phi_n\wedge\sigma_T$.
\end{proposition}
Let $c:S^1_s\sr S^1_s\vee S^1_s$ be the map in $H_{\bullet}$
corresponding to the usual codiagonal on the (simplicial) circle. 
\begin{proposition}
\llabel{add1}
Let $F_{\bullet}$ be a pointed simplicial sheaf. Then for any $p,q$
the map
$$\tilde{H}^{p,q}(S^1_s\wedge F_{\bullet})\oplus \tilde{H}^{p,q}(S^1_s\wedge
F_{\bullet})\sr \tilde{H}^{p,q}(S^1_s\wedge F_{\bullet})$$
defined by the codiagonal $c$ is of the form $(a,b)\mapsto a+b$.
\end{proposition}
\begin{proof}
It follows by adjunction argument from the fact that the map of the
free abelian groups defined by $c$ is isomorphic, in the derived
category of sheaves of abelian groups, to the diagonal map $\zz[1]\sr
\zz[1]\oplus\zz[1]$.
\end{proof}
\begin{cor}
\llabel{add2}
Let $\alpha:\tilde{H}^{p,q}\sr \tilde{H}^{r,s}$ be a cohomological operation. Then for
any pointed simplicial sheaf $F_{\bullet}$ the map:
$$\tilde{H}^{p,q}(S^1_s\wedge F_{\bullet})\sr \tilde{H}^{r,s}(S^1_s\wedge
F_{\bullet})$$
defined by $\alpha$ is a homomorphism of abelian groups.
\end{cor}
\begin{proof}
Follows from Proposition \ref{add1} using the naturality of $\alpha$
with respect to the map defined by the codiagonal $c$.
\end{proof}
\begin{cor}
\llabel{add3}
Let $\alpha:\tilde{H}^{*,*}\sr \tilde{H}^{*+i,*+j}$ be a bistable cohomological
operation. Then for any $F_{\bullet}$ the map
$$\tilde{H}^{*,*}(F_{\bullet})\sr \tilde{H}^{*+i,*+j}(F_{\bullet})$$
defined by $\alpha$ is a homomorphism of abelian groups.
\end{cor}
\begin{proof}
Follows from Corollary \ref{add2}.
\end{proof}

\subsection{Operations $\tilde{H}^{2d,d}\sr \tilde{H}^{2d+*,d+i}$ for $i\le 0$}\label{sec3}
In this section it will be convenient for us to use a different model
for the space $K_{n,A}$. We define $K_{n,A}'$ as the sheaf which sends
$U$ to the group of cycles with coefficients in $R$ on $U\times{\bf
A}^n$ which are equidimensional and of relative dimension zero over
$U$. The following theorem is proved in \cite{comparison}.
\begin{theorem}
\llabel{comp} There is an isomorphism $K_{n,A}\sr K_{n,A}'$ in the
$\af$-homotopy category of sheaves with transfers.
\end{theorem}
As a corollary we get the following result.
\begin{cor}
\llabel{comp2}
The pointed sheaves $K_{n,A}$ and $K_{n,A}'$ are isomorphic in
$H_{\bullet}$.
\end{cor}
Given a pointed sheaf $F$ define its standard simplicial resolution as
the simplicial sheaf $G_{\bullet}F$ with terms of the form
$$G_{i}F=(\coprod_{X_0\sr\dots\sr X_i; f\in F(X_i)-\{*\}}X_0)_+$$
where the coproduct is taken over all sequences of morphisms of length
$i$ in some small subcategory equivalent to $Sm/k$ (see
\cite{delnotes} for more details). One verifies easily that the
obvious morphism $G_{\bullet}F\sr F$ is a weak equivalence of pointed
simplicial presheaves. 

For a cycle $Z$ on $X$ denote by $Supp(Z)$ the closure of the set of
points which appear in $Z$ with nonzero multiplicity.  Consider
$G_{\bullet}K'_{n,A}\times{\bf A}^n$. For each $i$ let $F_{n,i}$ be
the open subset in
$$G_{i}K'_{n,A}\times{\bf A}^n=(\coprod_{X_0\sr\dots\sr X_i; Z\in
z(X_i\times{\bf A}^n/X_i)-\{0\}}{X_0\times {\bf A}^n})_+$$
whose component corresponding to $(X_0\sr\dots\sr X_i; Z)$ is the
complement to $X_0\times_{X_i} Supp(Z)$. The following lemma is
straightforward.
\begin{lemma}
The collection of subsheaves $F_{n,i}$ forms a simplicial subsheaf in
$G_{\bullet}K'_{n,A}\times{\bf A}^n$. 
\end{lemma}
\begin{proposition}
\llabel{hom}
The composition $F_{n,\bullet}\sr G_{\bullet}K'_{n,A}\times{\bf
A}^n\sr K'_{n,A}$ is $\af$-homotopic to the zero morphism.
\end{proposition}
\begin{proof}
To prove the proposition it is sufficient to construct for any $X$ and
any $Z$ in $Hom(X,K'_{n,A})=z(X\times{\bf A}^n/X)$ an $\af$-homotopy
from the map
$$X\times {\bf A}^n-Supp(Z)\sr X\sr K'_{n,A}$$
such that these homotopies are natural in $X$. Consider the map
$$h:(X\times {\bf A}^n-Supp(Z))\times {\bf A}^n\times\af\sr (X\times
{\bf A}^n-Supp(Z))\times {\bf A}^n$$
which sends $(x,u,v,t)$ to $(x,u,u(1-t)+vt$. This map is flat over the
complement to $X\times \Delta({\bf A}^n)$ in the target. Consider the
pull-back $p^*(Z)$ of $Z$ along the map $p:X\times {\bf
A}^n-Supp(Z)\sr X$. It is a cycle on the target of $h$ and the support
of this cycle does not intersect $X\times\Delta({\bf A}^n)$.
Therefore, the flat pull-back $h^*(p^*(Z))$ is defined. One verifies
easily that this cycle is equidimensional over $(X\times {\bf
A}^n-Supp(Z))\times\af$ of relative dimension zero and hence defines a
map:
$$H:(X\times {\bf A}^n-Supp(Z))\times\af\sr K'_{n,A}$$
The restriction of $h^*(p^*(Z))$ to $t=0$ is $p^*(Z)$ and the
restriction to $t=1$ is the zero cycle. Therefore, $H$ is a homotopy
of the required form. It is clear that our construction is
functorial in $X$.
\end{proof}
\begin{cor}
\llabel{retr}
The object $K_{n,A}'$ is a retract, in $H_{\bullet}$, of
$G_{\bullet}K'_{n,A}\times{\bf A}^n/F_{n,\bullet}$.
\end{cor}
We will often use below the following result.
\begin{lemma}
\llabel{often}
Let $k$ be a field, $X$ be smooth scheme over $k$ and $Z$ a
closed subscheme in $X$ everywhere of codimension at least $c$. Then
$\tilde{H}^{*,q}(X/(X-Z))=0$ for $q<c$ and
$$
\tilde{H}^{p,c}(X/(X-Z),A)=\left\{
\begin{array}{cl}
\oplus_{z\in Z^{c}}A &\mbox{\rm for $p=2c$}\\
0 &\mbox{\rm for $p\ne 2c$}
\end{array}
\right.
$$
where $Z^c$ is the set of points of $Z$ which are of codimension $c$
in $X$.
\end{lemma}
\begin{proposition}
\llabel{neg}
For any $n\ge 0$, $m<n$ and any abelian groups $A, B$ one has
$$\tilde{H}^{*,m}(K_{n,A}, B)=0$$
\end{proposition}
\begin{proof}
Follows immediately from Corollary \ref{retr} and Lemma \ref{often}.
\end{proof}
\begin{proposition}
\llabel{zer}
For any $n>0$ and any abelian groups $A, B$ one has
$$\tilde{H}^{p,n}(K_{n,A}, B)=\left\{
\begin{array}{cl}
Hom(A,B) &\mbox{\rm for $p=2n$}\\
0 &\mbox{\rm for $p<2n$}
\end{array}
\right.
$$
\end{proposition}
\begin{proof}
The fact that $\tilde{H}^{p,n}(K_{n,A}, B)=0$ for $p<2n$ follows immediately
from Corollary \ref{retr} and Lemma \ref{often}. Consider the case
$p=2n$.  We have an obvious map from $Hom(A,B)$ to $\tilde{H}^{2n,n}(K_{n,A},
B)$. On the other hand, an element of $\tilde{H}^{2n,n}(K_{n,A},
B)$ considered as an operation defines a map 
$$A=\tilde{H}^{2n,n}(T^n,A)\sr \tilde{H}^{2n,n}(T^n,B)=B$$
Since $T=S^1_s\wedge S^1_t$ and $n>0$, Corollary \ref{add2} implies
that this map is a homomorphism. Therefore, it is sufficient to show
that an element $\alpha$ of $\tilde{H}^{2n,n}(K_{n,A}, B)$ which acts
trivially on $\tilde{H}^{2n,n}(T^n,A)$ is zero. By Corollary \ref{retr} any
$\alpha$ is determined by its action on objects of the form
$X\times{\bf A}^n/(X\times{\bf A}^n-Z)$ where $Z$ is a closed subset
equidimensional of relative dimension $0$ over ${\bf A}^n$. Let
$Z_{sing}$ be the closed subset of singular points of $Z$. Since $k$
is perfect, $Z_{sing}$ is of codimension at least $n+1$ in
$X\times{\bf A}^n$. Lemma \ref{often} implies that the motivic
cohomology of weight $n$ of $X\times{\bf A}^n/(X\times{\bf A}^n-Z)$
map to the motivic cohomology of weight $n$ of $(X\times{\bf
A}^n-Z_{sing})/(X\times{\bf A}^n-Z)$ isomorphically. It remains to
show that $\alpha$ acts trivially on
$$\tilde{H}^{2n,n}((X\times{\bf
A}^n-Z_{sing})/(X\times{\bf A}^n-Z),A)$$
The normal bundle to $Z-Z_{sing}$ in $X\times{\bf
A}^n-Z_{sing}$ is trivial. Hence, by the homotopy purity theorem
\cite{MoVo}, we have a weak equivalence
$$(X\times{\bf
A}^n-Z_{sing})/(X\times{\bf A}^n-Z)=\Sigma^n_T((Z-Z_{sing})_+)$$
let $Z_i$, $i=1,\dots, m$ be the connected components of
$Z-Z_{sing}$. Then we have a map 
$$\Sigma^n_T((Z-Z_{sing})_+)\sr \vee_{i=1}^m T^n$$
and $H^{j,0}(Z_i)$ is non-zero only for $j=0$ where it is $A$, this
map defines an isomorphism on $\tilde{H}^{2n,n}(-,A)$. We conclude that
$\alpha$ acts by zero since by assumption it acts by zero on
$\tilde{H}^{2n,n}(T^n,A)$.  
\end{proof}

\subsection{Thom isomorphism and Euler classes}\llabel{sec4}
If $E$ is a vector bundle and ${\bf P}(E)$ is the projective bundle
defined by $E$ then the line bundle ${\cal O}(-1)$ on ${\bf P}(E)$
gives a class in $H^{2,1}({\bf P}(E),\zz)$ which we denote by
$\sigma$. The following result is proved in \cite{}[].
\begin{theorem}
\llabel{projbund}\llabel{l20}
For any smooth $X$ over $k$ and a vector bundle $E$ on $X$ of
dimension $d$, the elements $1,\sigma,\dots,\sigma^{d-1}$ form a
basis of the $H^{*,*}$-module $H^{*,*}({\bf P}(E),\zz)$.
\end{theorem} 
The key ingredient of the proof of this theorem is the following lemma
which we will also use directly.
\begin{lemma}
\llabel{morel}
Let $0$ be the image of the point $(0,\dots,0)$ under the standard
embedding ${\bf A}^n\sr {\bf P}^n$. Let further $f:{\bf P}^n\sr T^n$
be the composition 
$${\bf P}^n\sr {\bf P}^n/({\bf P}^n-0)\cong {\bf A}^n/({\bf
A}^n-0)=T^n$$
and $t$ the tautological class in $\tilde{H}^{2n,n}(T^n,\zz)$. Then
$f^*(t)=(-\sigma)^n$.
\end{lemma}
Recall that for a vector bundle $E$ on $X$ we denote by $Th(E)$ the
pointed sheaf $E/(E-z(X))$ where $z:X\sr E$ is the zero
section. Consider the projective bundle ${\bf P}(E\oplus {\cal
O})$. We have two morphisms
$$X={\bf P}({\cal O})\sr {\bf P}(E\oplus {\cal
O})$$
$${\bf P}(E)\sr {\bf P}(E\oplus {\cal
O})$$
The complement to the image of the second morphism is $E$ and open
excision implies that 
$$Th(E)={\bf P}(E\oplus {\cal
O})/{\bf P}(E\oplus {\cal
O})-X$$
On the other hand the map ${\bf P}(E)\sr {\bf P}(E\oplus {\cal
O})-X$ is locally of the form ${\bf P}^{n-1}\sr {\bf P}^{n}-pt$ and
therefore it is an $\af$-weak equivalence. We conclude that the
morphism
$${\bf P}(E\oplus {\cal
O})/{\bf P}(E)\sr Th(E)$$
is a weak equivalence. If $f:E\sr E'$ is a monomorphism of vector
bundles and ${\bf P}(f):{\bf P}(E)\sr {\bf P}(E')$ is the
corresponding morphism of projective bundles then ${\bf P}(f)^*({\cal
O}(-1))={\cal O}(-1)$. Together with Theorem \ref{projbund} this
implies that the map on motivic cohomology defined by ${\bf P}(E)\sr
{\bf P}(E\oplus {\cal O})$ is a split mono and that there is a unique
class in $\tilde{H}^{2d,d}(Th(E),\zz)$ whose image in the cohomology
of ${\bf P}(E\oplus {\cal O})$ is of the form
$(-\sigma)^{d}+\sum_{i<d} a_i \sigma^i$. This class is called
the Thom class of $E$ and denoted $t_E$. 

The obvious ``diagonal'' map $d:Th(E)\sr X_+\wedge Th(E)$ defines
multiplication $(x,y)\mapsto d^*(x\wedge y)$ with $x\in H^{*,*}(X)$,
$y\in \tilde{H}^{*,*}(Th(E))$ and $d^*(x\wedge y)\in
\tilde{H}^{*,*}(Th(E))$. By abuse of notation we will write $x
y$ instead of $d^*(x\wedge y)$.
\begin{proposition}
\llabel{thom}\llabel{l7}
For any pointed simplicial sheaf $F_{\bullet}$ the map $a\mapsto
a t_E$ from $\tilde{H}^{*,*}(F_{\bullet}\wedge X_+)$ to
$\tilde{H}^{*+2d,*+d}(F_{\bullet}\wedge Th(E))$ is an isomorphism.
\end{proposition}
\begin{proof}
A standard argument shows that it is sufficient to prove the
proposition for $F_{\bullet}=pt$. In this case it follows immediately
from out definition of the Thom class and the projective bundle
theorem.
\end{proof}
\begin{cor}
\llabel{uniq}\llabel{l39}
The Thom class $t_E$ is a unique class in $\tilde{H}^{2d,d}(Th(E))$
whose restriction to any generic point of $X$ is the tautological
class in $\tilde{H}^{2d,d}(T^{d})$.
\end{cor}
\begin{proof}
The fact that the restriction of $t_E$ to the generic point is the
tautological class follows from Lemma \ref{morel}. The fact that $t_E$
is determined by this condition follows from Proposition \ref{thom}
and the fact that if $j:\coprod Spec(K_i)\sr X$ is the embedding of
the generic points of $X$, then $j^*$ defines an isomorphism on
$H^{0,0}$.
\end{proof}
For a vector bundle $E$ define the Euler class $e(E)$ in
$H^{2d,d}(X)$ as the restriction of $t_E$ with respect to
the zero section map $X_+\sr Th(E)$.
\begin{lemma}
\llabel{line}\llabel{l54}
Let $L$ be a line bundle. Then $e(L)$ coincides with the canonical
class of $L$ in $H^{2,1}$. In particular, for two line bundles $L$,
$L'$ one has $e(L\oo L')=e(L)+e(L')$.
\end{lemma}
\begin{proof}
Let $\sigma+c$ be the image of $t_L$ in $H^{2,1}({\bf P}(L\oplus {\cal
O}))$. The restriction of $\sigma$ to ${\bf P}({\cal O})$ is ${\cal
O}(-1)$ for ${\cal O}$ i.e. zero. Therefore, $e(L)=c$. On the
other hand $c$ is defined by the condition that the restriction of
$\sigma+c$ to ${\bf P}(L)$ is zero. Since the restriction of $\sigma$
is the class of $-L$ we conclude that $e(L)=L$. 
\end{proof}
\begin{lemma}
\llabel{sum}\llabel{l29}
Let $E, E'$ be two vector bundles. Then $e(E\oplus E')=e(E)e(E')$.
\end{lemma}
\begin{proof}
Corollary \ref{l39} implies that the Thom class for the sum of two
bundles is the smash product of Thom classes. In particular,
$e(E\oplus E')=e(E)e(E')$.
\end{proof}
\begin{lemma}
\llabel{mono}\llabel{l14} Let $f:E\sr E'$ be a monomorphism of vector
bundles on a quasi-projective scheme $X$ such that $E'/E$ is again a
vector bundle and $th(f):Th(E)\sr Th(E')$ be the corresponding map of
Thom spaces. Then $th(f)^*(t_E')=t_E e(E/E')$.
\end{lemma}
\begin{proof}
Since $X$ is quasi-projective we can find an affine torsor $X'\sr X$
such that $X'$ is affine. Since motivic cohomology of $X$ and $X'$ are
the same it is sufficients to prove the lemma for an affine $X$. Over
an affine $X$ the sequence $E\sr E'\sr E'/E$ splits and we have
$E'=E\oplus E'/E$. Let $d:Th(E\oplus E'/E)\sr Th(E)\wedge Th(E'/E)$ be
the obvious morphism. Corollary \ref{l39} implies that
$t_{E'}=d^*(t_E\wedge t_{E'/E})$ which in turn implies the statement
of the lemma.
\end{proof}

\subsection{Total power operations}
\begin{constr}\rm
\llabel{c1}
Let $E, L$ be vector bundles on $X$ and $\phi:E\oplus L\sr {\cal O}^n$
an isomorphism. For an cycle $Z$ on $E$ with coefficients in a
commutative ring $R$ which is equidimensional and finite over $X$
consider the cycle on $L\times_X E\times_X L$ whose fiber over a point
$(x,l)$ of $L$ is $Z_x\times l$. One verifies easily that this is a
cycle equidimensional and finite over $L$. Identifying $E\times_X L$
with ${\bf A}^n_X$ by means of $\phi$ we get a section of
$\zz_{tr}({\bf A}^N)$ on $L$. The restriction of this section to
$L-z(X)$ where $z:X\sr L$ is the zero section lies in $\zz_{tr}({\bf
A}^N-\{0\})$. Therefore, it gives a map of pointed sheaves $Th(L)\sr
K_{n,R}$ which we denote $a(Z)$. 
\end{constr}
\begin{lemma}
\llabel{p2}
The motivic cohomology class $\tilde{a}(Z)$ in $H^{2dim(E),dim(E)}(X)$
defined by $a(Z)$ through the Thom isomorphism does not depend on the
choice of $L$ and $\phi$.
\end{lemma}
\begin{proof}
Let $x$ be a generic point of $X$, $L_x$ the fiber of $L$ over $x$ and
$\sum z_i$ the fiber of $Z$ over $x$. Then, for a generic point $l$ of
$L_x$ the fiber of $\tilde{a}(Z)$ over $(x,l)$ is $\sum
\phi(z_i,l)$. Given another pair $(L',\phi')$ consider the isomorphism
$$L\oplus {\cal O}^{N'}\sr L\oplus E\oplus L'\sr L'\oplus E\oplus L\sr
L'\oplus {\cal O}^N$$
and let $\psi:Th(L\oplus {\cal O}^{N'})\sr Th(L\oplus {\cal O}^{N})$
be the corresponding isomorphism of Thom spaces. We claim that
\begin{equation}
\llabel{oncycles}
\psi^*(a(Z) t_N')=a'(Z) t_{N}
\end{equation}
on the level of actual cycles. Indeed, the fiber of $a(Z)
t_{N'}$ over a generic point $(x,l,u)$ of $L\oplus {\cal O}^{N'}$ is
$\sum (\phi(z_i,l),u)$ (in ${\cal O}^N\oplus {\cal O}^{N'}$) which
coincides with the fiber of $a'(Z) t_{N}$ over
$\psi(x,l,u)$. Using (\ref{oncycles}) and applying Lemma \ref{mono} to
$\psi^*$ we get the statement of our lemma.
\end{proof}
\begin{constr}\rm
\llabel{c2}
Let $G$ be a finite group, $r:G\sr S_n$ a permutational
representation of $G$, $U$ a smooth scheme with free action of $G$ and
$L$ a vector bundle on $U/G$ given together with an isomorphism
$\xi_n\oplus L\sr {\cal O}^N$ where $\xi_n$ is the vector bundle of
dimension $n$ on $U/G$ corresponding to $r$. Given any such collection
and a cycle $Z$ on $X\times{\bf A}^i$ equidimensional and finite over
$X$ define a map
$$\tilde{P}(Z):X\wedge Th_{U/G}(L^i)\sr K_{iN,R}$$
as follows. Let $Z^{\oo n}$ be the external power of $Z$. It is a
cycle on $(X\times{\bf A}^i)^n$. Let $p^*(Z^{\oo n})$ be its flat
pull-back to $(X\times{\bf A}^i)^n\times U$. Since $p^*(Z^{\oo n})$ is
invariant under the action of $G$ and the action of $G$ on $U$ is free
there exists a unique cycle $Z'$ on $((X\times{\bf A}^i)^n\times U)/G$
whose pull-back to $(X\times{\bf A}^i)^n\times U$ is $Z$. One verifies
easily that $Z'$ is finite and equidimensional over $(X^n\times
U)/G$. Therefore, we can pull it back to a cycle $Z''$ on $X\times
({\bf A}^{in}\times U)/G$ by means of the diagonal map $X\sr X^i$. The
scheme $({\bf A}^{in}\times U)/G$ is the vector bundle $\xi_n^i$ over
$U/G$ and we define $\tilde{P}(Z)$ as $a(Z'')$. One verifies
immediately that if $Z$ lies in $X\times({\bf A}^i-\{0\})$ then
$\tilde{P}(Z)=0$. 
\end{constr}
This construction defines a morphism of pointed sheaves:
$$\tilde{P}:K_{i,R}\wedge Th_{U/G}(L^i)\sr K_{iN,R}$$
Since smash products preserve $\af$-weak equivalences, this morphism
defines an operation
$$\tilde{P}:\tilde{H}^{2i,i}(-,R)\sr \tilde{H}^{2iN,iN}(-\wedge Th_{U/G}(L^i),R)$$
By Lemma \ref{l7} the Thom isomorphism defines a morphism in
$H_{\bullet}$ of the form
$$P:K_{i,R}\wedge(U/G)_+\sr K_{in,R}$$
or, equivalently an operation
$${P}:\tilde{H}^{2i,i}(-,R)\sr \tilde{H}^{2in,in}(-\wedge (U/G)_{+},R)$$
such that 
\begin{equation}
\llabel{use}
\tilde{P}(x)=P(x) t_{L^i}
\end{equation}
Lemma \ref{p2} implies
immediately the following result.
\begin{lemma}
\llabel{l30}
The operation $P=P_{G,r,U,L,\phi}$ does not depend on the choice of
$L$ and $\phi$.
\end{lemma}
If $U$ is a quasi-projective scheme with a free action of $G$ then we
can find, using the standard trick, an affine smooth scheme
$\tilde{U}$ with a free action of $G$ and an equivariant morphism
$\tilde{U}\sr U$. Since any vector bundle $E$ on an affine scheme is
``invertible'', i.e. there is an $L$ such that $L\oplus E\cong {\cal
O}^N$, we can define $P$ for $\tilde{U}$. Since $\tilde{U}\sr U$ is an
$\af$-weak equivalence this means that we have a well defined
operation $P$ for any $G$, $r:G\sr S_n$ and any quasi-projective $U$
with a free $G$-action. The following lemma is straightforward.
\begin{lemma}
\llabel{l33}
Let $G$ be a finite group and $r:G\sr S_n$ a permutational
representation of $G$. Let $U$, $V$ be smooth quasi-projective schemes
with free actions of $G$ and $f:U\sr V$ an equivariant morphism. Then
for any $x$ one has $P_{U}(x)=f^*(P_{V}(x))$.
\end{lemma}
\begin{lemma}
\llabel{l21}
Let $G, r, U, L, \phi$ be as above. Then the following diagram of
morphisms of pointed sheaves commutes
\begin{equation}
\llabel{prod}
\begin{CD}
K_{i}\wedge K_j\wedge Th_{U/G}(L^{i+j}) @>>> K_{i}\wedge
Th_{U/G}(L^{i})\wedge K_j\wedge
Th_{U/G}(L^{j})\\
@VVV @VVV\\
@VVV K_{iN}\wedge K_{jN}\\
@VVV @VVV\\
K_{i+j}\wedge Th_{U/G}(L^{i+j}) @>>> K_{(i+j)N}
\end{CD}
\end{equation}
\end{lemma}
\begin{proof}
Direct comparison.
\end{proof}
\begin{lemma}
\llabel{l22}
For $a\in \tilde{H}^{2i,i}(F_{\bullet})$ and $b\in \tilde{H}^{2j,j}(F'_{\bullet})$ one
has 
\begin{equation}
\llabel{prodf}
P(a\wedge b)=\Delta^*(P(a)\wedge P(b))
\end{equation}
where $\Delta:U/G\sr U/G\times U/G$ is the diagonal.
\end{lemma}
\begin{proof}
By Lemma \ref{prod} we have $\tilde{P}(a\wedge
b)=\delta^*(\tilde{P}(a)\wedge \tilde{P}(b))$ where 
$$\delta:Th_{U/G}(L^{i+j}) \sr  Th_{U/G}(L^{i})\wedge
Th_{U/G}(L^{j})$$
Therefore by (\ref{use}):
$${P}(a\wedge b) t_{L^{i+j}}=\delta^*(({P}(a) t_{L^i})\wedge
({P}(b) t_{L^j}))$$
Since $\delta^{*}(t_{L^i}\wedge t_{L^j})=t_{L^{i+j}}$ the Thom
isomorphism theorem implies (\ref{prodf}).
\end{proof}
\begin{lemma}
\llabel{l40}
Let $\eta^i\in \tilde{H}^{2i,i}(T^i)$ be the tautological class. Then
$P(\eta^i)=\delta^*(t_{\xi_n^i})$ where $\xi_n$ is the vector bundle
on $U/G$ corresponding to the representation $r:G\sr S_n$ and $\delta$
is the map on Thom spaces defined by the embedding of vector bundles
${\cal O}^i\sr \xi_n^i$ on $U/G$.
\end{lemma}
\begin{proof}
Take ${\bf A}^i/{\bf A}^i-\{0\}$ as a model of $T^i$ such that
$\eta^i$ is given by the tautological section $Z$ of $\zz_{tr}({\bf
A}^i)$ on ${\bf A}^i$. Applying the construction of $\tilde{P}$ to the
corresponding diagonal cycle $Z$ we get the restriction to $Th(L^i)$
of the morphism
$$W:Th(L^i\oplus \xi^i_n)\sr \zz_{tr}({\bf A}^Ni)$$
corresponding, by Construction \ref{c1}, to the tautological cycle on
$\xi_n^{2i}$ over $\xi^i_n$. In view of Corollary
\ref{l39} $W$ represents the Thom class of the trivial bundle
$(L\oplus\xi_n)^i$. Applying Thom isomorphism to get $P(Z)$ we
conclude that $P(Z)$ is the restriction of the Thom class of
$Th(\xi^i_n)$ with respect to the morphism $\delta$.
\end{proof}
\begin{lemma}
\llabel{l8}
Let $E$ be a vector bundle on $X$. Then
$P(t_E)=\delta^*(t_{E\oo\xi_n})$ where $\xi_n$ is the vector bundle on
$U/G$ corresponding to the representation $r:G\sr S_n$ and $\delta$ is
the map on Thom spaces defined by the embedding of vector bundles
$E\oo {\cal O}\sr E\oo \xi_n$ on $X\times(U/G)$.
\end{lemma}
\begin{proof}
Follows from Corollary \ref{l39} and Lemma \ref{l40}.
\end{proof}
Let $*$ be a $k$-point of $U/G$ which lifts to a $k$-point of $U$. Let
$i:S^0\sr (U/G)_+$ be the corresponding morphism. The following lemma
is straightforward.
\begin{lemma}
\llabel{power0}
The composition
$$K_{i,R}\stackrel{Id\wedge i}{\sr}K_{i,R}\wedge
(U/G)_+\stackrel{P}{\sr} K_{in,R}$$
coincides in $H_{\bullet}$ with the n-th power map.
\end{lemma}
\begin{proof}
Let $L_{*}$ be the fiber of $L$ over our distinguished point
$*$. Then, if we compute the analog of our composition using
$\tilde{P}$ instead of $P$, we get the map $K_{i,R}\wedge Th(L_*)\sr
K_{iN,R}$ which is of the form $a(Z^{\oo n})$ where $Z$ is the
tautological cycle on $K_{i,R}$ and $a(Z)$ is Construction \ref{c1}
with respect to the isomorphism $(\xi_n)^{\oplus i}_*\oplus L^{\oplus
i}_*\sr {\cal O}^{iN}$. Our result follows now from Lemma \ref{p2}
since $(\xi_n)_*={\cal O}^{n}$ is the trivial bundle.
\end{proof}

\subsection{Motivic cohomology of $B\mu_l$ and $BS_l$}\label{sec6}
Let $G$ be a linear algebraic group and $G\sr GL(V)$ a faithful
representation of $G$. Denote by $\tilde{V}_n$ the open subset in ${\bf
A}(V)^n$ where $G$ acts freely. We have a sequence of closed
embeddings $\tilde{V}_n\sr \tilde{V}_{n+1}$ given by
$(v_1,\dots,v_n)\mapsto (v_1,\dots,v_n,0)$. Set $BG=colim_n
\tilde{V}_n/G$ where $\tilde{V}_n/G$ is the quotient scheme and the
colimit is taken in the category of sheaves. In \cite{MoVo} we used
the notation $B_{gm}G$ for $BG$. As shown there, the homotopy type of
$BG$ does not depend on the choice of $G\sr GL(V)$. We denote by $*$
any $k$-rational point of $BG$ which lifts to a $k$-rational point in one of
the $\tilde{V}_n$'s. The goal of this section is to describe motivic
cohomology of $BS_l$ with coefficients in $\zz/l$. We start with
the following general result.
\begin{proposition}
\llabel{newprop}
For any $G$ and $V$ as above the map 
$$i_n:\tilde{V}_{n}/G\sr
\tilde{V}_{n+1}/G$$
defines an isomorphism on motivic cohomology of weigh less than $n$.
\end{proposition}
\begin{proof}
The morphism $\tilde{V}_{n}/G\sr
\tilde{V}_{n+1}/G$ can be represented by the composition
$$\tilde{V}_n/G\sr (\tilde{V}_n\times V)/G\sr \tilde{V}_{n+1}/G$$
The first of these maps is the zero section of a vector bundle
$(\tilde{V}_n\times V)/G\sr\tilde{V}_n/G$ and gives an isomorphism on
motivic cohomology by homotopy invariance.  The second map is an open
embedding of smooth schemes and the codimension of the complement is
at least $n$. By Lemma \ref{often} it defines an isomorphism on
$\tilde{H}^{*,<n}$.
\end{proof}
\begin{cor}
\llabel{limit}
One has
$$H^{*,*}(BG)=lim_{n} H^{*,*}(\tilde{V}_m/G)$$
\end{cor}
\begin{proof}
Given any sequence of maps of pointed simplicial sheaves
$F_{\bullet,n}\sr F_{\bullet,n+1}$ with the colimit
$F_{\bullet,\infty}$ we have a long exact sequence of the form
$$\sr\prod \tilde{H}^{*-1,*}(F_{\bullet,n}) \sr
\tilde{H}^{*,*}(F_{\bullet,\infty})\sr \prod \tilde{H}^{*,*}(F_{\bullet,n})\sr \prod
\tilde{H}^{*,*}(F_{\bullet,n})\sr$$
The limit $lim_n \tilde{H}^{*,*}(F_{\bullet,n})$ is the kernel of the fourth
arrow and therefore to prove the corollary it is sufficient to show that 
the map 
$$Id-\prod i_n^*:\prod H^{*,*}(\tilde{V}_n/G)\sr \prod
H^{*,*}(\tilde{V}_n/G)$$
is an epimorphism. It follows from the proposition.
\end{proof}
Denote by $\mu_l$ the groups scheme of l-th roots of unity
$$\mu_l:=ker({\bf G}_m\stackrel{z^l}{\sr}{\bf G}_m)$$
\begin{lemma}
\llabel{l43}
Let $B\mu_l$ be defined with respect to the tautological 1-dimensional
representation of $\mu_l$. Then one has
\begin{equation}
\llabel{iso}
B{\bf \mu}_l={\cal O}(-l)_{{\bf P}^{\infty}}-z({\bf P}^{\infty})
\end{equation}
\end{lemma}
\begin{proof}
We have $\tilde{V}_n={\bf A}^n-\{0\}$. The projection ${\bf
A}^n-\{0\}\sr {\bf P}^{n-1}$ is invariant under the action of $\mu_l$
and therefore gives a map $({\bf A}^n-\{0\})/\mu_l\sr {\bf
P}^{n-1}$. One verifies that this map is isomorphic to the complement
to the zero section of the line bundle ${\cal O}(-l)$ on ${\bf
P}^{n-1}$.
\end{proof}
Let $r:G\sr GL(V)$ be a linear representation of $G$. It defines an
action of $G$ on the affine space 
$${\bf A}(V)=Spec(S^{\BB}V^*)$$
corresponding to $V$. If $U$ is a scheme with a free action of $G$
then the projection $({\bf A}(V)\times U)/G\sr U/G$ is a vector
bundle. We say that this is the vector bundle defined by $r$.
\begin{lemma}
\llabel{corres} The line bundle on $B\mu_l$ defined by the
tautological representation of $\mu_l$ is isomorphic
with respect to (\ref{iso}) to the pull-back of ${\cal O}(1)$.
\end{lemma}
\begin{proof}
Let $L=\af$ with the standard action of ${\bf G}_m$ and the
corresponding action of $\mu_l$. The square
$$
\begin{CD}
(L\times({\bf A}^n-\{0\}))/\mu_l @>>> (L\times({\bf A}^n-\{0\}))/{\bf
G}_m\\
@VVV @VVV\\
({\bf A}^n-\{0\})/\mu_l @>>> ({\bf A}^n-\{0\})/{\bf
G}_m
\end{CD}
$$
is pull-back. The fact that the right vertical arrow is ${\cal
O}(1)\sr {\bf P}^{n-1}$ is standard (e.g. it is not ${\cal O}(-1)$
because it has a section other than the zero one).
\end{proof}
Lemma \ref{l43} implies that one has a cofibration sequence of the
form
\begin{equation}
\llabel{cs1}
(B{\bf \mu}_l)_+\sr ({\cal O}(-l)_{{\bf P}^{\infty}})_+\sr Th({\cal
O}(-l))
\end{equation}
For a vector bundle $E$ of dimension $d$, the composition of the Thom
isomorphism $H^{*,*}(X)\sr \tilde{H}^{*+2d,*+d}(Th(E))$ with the restriction to the
zero section $z^*:\tilde{H}^{*,*}(Th(E))\sr H^{*,*}(X)$ is given by 
$$x\mapsto
z^*(x t_E)=x z^*(t_E)=x e(E)$$
By Lemma \ref{l54}, $e({\cal O}(-l))=l\sigma$ where $\sigma\in
H^{2,1}({\bf P}^{\infty})$ is the same class as in the projective
bundle theorem \ref{l20}. Therefore, the long exact sequence defined
by (\ref{cs1}) is of the form
\begin{equation}
\llabel{mainseq}
\dots\sr H^{*-2,*-1}[[\sigma]]\stackrel{
l\sigma}{\sr}H^{*,*}[[\sigma]]\sr H^{*,*}(B{\bf \mu}_l)\sr
H^{*-1,*-1}[[\sigma]]\sr\dots
\end{equation}
The short exact sequence of abelian groups
$$0\sr \zz\sr\zz\sr\zz/l\sr 0$$
defines a homomorphism
$$\delta:\tilde{H}^{*,*}(-,\zz/l)\sr \tilde{H}^{*+1,*}(-,\zz)$$
Let $v$ be Euler class of the line bundle on $B\mu_l$ corresponding to
the tautological representation of $\mu_l$.
\begin{lemma}
\llabel{defu} There exists a unique element $u\in H^{1,1}(B{\bf
\mu}_l,\zz/l)$ such that the restriction of $u$ to $*$ is zero and
$\delta(u)=v$.
\end{lemma}
\begin{proof}
Existence follows from the exact sequence
$$H^{1,1}(B{\bf \mu}_l,\zz)\sr H^{1,1}(B{\bf \mu}_l,\zz/l)\sr
H^{2,1}(B{\bf \mu}_l,\zz)\stackrel{l}{\sr}H^{2,1}(B{\bf \mu}_l,\zz)$$
and the fact that 
$$lv=-le({\cal O}(-1))=-e({\cal O}(-l))=0$$
in $H^{2,1}(B{\bf \mu}_l,\zz)$. The exact sequence (\ref{mainseq})
around $H^{1,1}(B{\bf \mu}_l,\zz)$ shows that $H^{1,1}(B{\bf
\mu}_l,\zz)=H^{1,1}(k)$. This implies the uniqueness.
\end{proof}
\begin{proposition}
\llabel{l46}
For any pointed simplicial sheaf $F_{\bullet}$ the elements $v^i$ and
$uv^i$, $i\ge 0$ form a basis of $\tilde{H}^{*,*}(F_{\bullet}\wedge(B{\bf
\mu}_l)_+)$ over $\tilde{H}^{*,*}(F_{\bullet})$.
\end{proposition}
\begin{proof}
The standard argument shows that it is sufficient to consider the case
when $F_{\bullet}$ is of the form $X_+$ for a smooth scheme $X$.  The
same reasoning as we used to establish (\ref{mainseq}) applies to the
motivic cohomology groups of $X\times B{\bf \mu}_l$ for any smooth
scheme $X$ and we get the following result.
\begin{lemma}
\llabel{new}
For any smooth scheme $X$ there is a long exact sequence of
$H^{*,*}(X)[[\sigma]]$-modules of the form
\begin{equation}
\llabel{main2}
\begin{array}{c}
\dots\sr H^{*-2,*-1}(X)[[\sigma]]\stackrel{
l\sigma}{\sr}H^{*,*}(X)[[\sigma]]\sr
H^{*,*}(X\times B{\bf \mu}_l)\sr\\\\
\sr H^{*-1,*-1}(X)[[\sigma]]\sr\dots
\end{array}
\end{equation}
\end{lemma}
For $\zz/l$-coefficients we have $l\sigma=0$ and (\ref{main2}) becomes
a short exact
sequence of $H^{*,*}(X)[[\sigma]]$-modules of the form 
\begin{equation}
\llabel{split}
0\sr H^{*,*}(X)[[\sigma]]\sr
H^{*,*}(X\times B{\bf \mu}_l)\sr
H^{*-1,*-1}(X)[[\sigma]]\sr 0
\end{equation}
Let $u'$ be an element in $H^{1,1}(B{\bf \mu}_l,\zz/l)$ such that the
image of $u'$ in 
$$H^{0,0}(Spec(k))=\zz/l$$
is $1$ and the restriction of $u'$ to $X\times *$ is zero. Since $v$
is the image of $\sigma$, the short exact sequence (\ref{split})
implies that the monomials $u'v^i$ and $v^i$ form a basis of
$H^{*,*}(X\times B{\bf \mu}_l)$ over $H^{*,*}(X)$. On the other hand,
the image of $u$ in $H^{0,0}(Spec(k))=\zz/l$ is not zero and hence
$u=cu'$ where $c\in (\zz/l)^*$. This implies that the monomials
$uv^i$, $v^i$ also form a basis.
\end{proof}
To describe the multiplicative structure of $\tilde{H}^{*,*}(F_{\bullet}\wedge(B{\bf
\mu}_l)_+)$ it is sufficient to find an expression for $u^2$. If $l\ne
2$ then $u^2=0$ since the multiplication in motivic cohomology is
graded commutative. Consider the case $l=2$. We can clearly assume
that $F_{\bullet}=Spec(k)_+$. The element $u^2$ lies in $H^{2,2}$ and
Lemma \ref{l46} shows that
$$H^{2,2}(B{\bf
\mu}_l)=H^{0,1}(k)v\oplus H^{1,1}(k)u\oplus H^{2,2}(k)$$
Since $u$ is zero in $*$ the projection of $u^2$ to the last factor is
zero and we get $u^2=xv+yu$ for $x\in H^{0,1}(k)$ and $y\in
H^{1,1}(k)$. To compute $y$ consider the map 
\begin{equation}
\llabel{manywords}
{\bf A}^1-\{0\}=({\bf A}^1-\{0\})/\mu_l\sr B{\bf \mu}_l
\end{equation}
If we choose the distinguished point of $B\mu_l$ to be the image of
the point $(1,0,\dots)$, this map is the embedding of the fiber of the
line bundle ${\cal O}(-l)-z({\bf P}^{\infty})\sr {\bf P}^{\infty}$
which contains the distinguished point. The pull-back along this map
coincides with the composition of the last map of (\ref{mainseq}) with
the map $H^{*,*}[[\sigma]]\sr \tilde{H}^{*,*}(\af-\{0\},1)$ which sends $1$ to
the generator of $\tilde{H}^{1,1}(\af-\{0\},1)$ and $\sigma$ to zero. In
particular, the pull-back of $u$ along (\ref{manywords}) is non-zero.
The following lemma implies now that $y=\rho$ where $\rho$ is the class of
$-1$ in $H^{1,1}$.
\begin{lemma}
\llabel{starn}\llabel{l48}
Let $w$ be the non zero element of 
$$\tilde{H}^{1,1}(({\bf
A}^1-\{0\},1),\zz/2)=\zz/2.$$
Then $w^2=\rho w$.
\end{lemma}
\begin{proof}
We need to compare two motivic cohomology classes in
$H^{2,2}(\af-\{0\})$. Let $Spec(k(t))\sr \af-\{0\}$ be the embedding
of the generic point. Since the base field $k$ may be assumed to be
perfect, the Gysin long exact sequence in motivic cohomology implies
that the kernel of the induced map in $H^{2,2}$ is covered by a direct
sum of groups of the form $H^{0,1}(Spec(E),\zz)$. Since such groups
are zero it is a monomorphism. Therefore, it is sufficient to show
that $t t=\rho  t$ in $H^{2,2}(Spec(k(t)),\zz)$. By
\cite{SusVoe3}, this group is isomorphic to $K_2^M(k(t))$ and we
conclude by the well known relation $(a,a)=(-1,a)$ in the Milnor's
K-theory. 
\end{proof}
To compute $x$ note that $H^{0,1}(k,\zz/2)={\bf \mu}_2(k)$. If
$char(k)=2$ then this group is zero. If $char(k)\ne 2$ it is $\zz/2$
and we only need to know whether $x$ is zero or not. The following
lemma implies that $x=\tau$ is the generator of $\zz/2$.
\begin{lemma}
\llabel{l49}
Let $k$ be a separably closed field of characteristic not equal to
$2$. Then $u^2\ne 0$.
\end{lemma}
\begin{proof}
We have a natural transformation from the motivic cohomology to the
etale cohomology with $\zz/2$-coefficients. For a class $u$ in the
etale $H^1$ we have $u^2=\beta(u)$ where $\beta$ is the Bockstein
homomorphism. Since $k$ is separably closed and in particular contains
$\sqrt{-1}$, the Bockstein in the etale cohomology commutes with the
Bockstein in the motivic cohomology and we conclude that the image of
$u^2$ in the etale cohomology coincides with the image of $v=\beta(u)$
in the etale cohomology. An etale analog of the long exact sequence
(\ref{mainseq}) shows that the image of $v$ in the etale cohomology is non
zero.
\end{proof}
We proved the following result.
\begin{theorem}
\llabel{l47}
For any field $k$ and a pointed simplicial sheaf $F_{\bullet}$ over
$k$ one has:
\begin{equation}
\llabel{eq47}
\tilde{H}^{*,*}(F_{\bullet}\wedge (B{\bf
\mu}_l)_+,\zz/l)=\tilde{H}^{*,*}(F_{\bullet},\zz/l)[[u,v]]/(u^2=\tau v +\rho u)
\end{equation}
where:
\begin{enumerate}
\item $\rho$ is the class of $-1$ in $H^{1,1}(k)$
\item $\tau$ is zero if $l\ne 2$ or $char(k)=2$
\item $\tau$ is the generator of $H^{0,1}(k,\zz/2)={\bf \mu}_2(k)$ if
$l=2$ and $char(k)\ne 2$.
\end{enumerate}
\end{theorem}
We also need the following additional fact about cohomology of $B{\bf
\mu}_l$. 
\begin{lemma}
\llabel{l52}
Let $c\in Aut({\bf \mu_l})=(\zz/l)^*$ and let $\uu{c}$ be the
corresponding automorphism of $B{\bf
\mu}_l$. Then $\uu{c}^*(u)=c u$ and $\uu{c}^*(v)=c v$.
\end{lemma}
\begin{proof}
Let $L$ be the line bundle on $B{\bf \mu}_l$ corresponding to the
tautological 1-dimensional representation $\lambda$ of $\mu_l$. Then
$v$ is the Euler class of $L$ and $u$ is the only element in
$H^{1,1}(-,\zz/l)$ which is zero at $*$ and which maps to $v$ under
the map $\delta:H^{1,1}(-,\zz/l)\sr H^{2,1}(-,\zz)$. The automorphism
$\uu{c}$ takes the tautological 1-dimensional representation $\lambda$
to $\lambda^{\oo c}$ and, therefore, it takes $L$ to $L^{\oo c}$. Our
result follows now from Lemma \ref{l54}.
\end{proof}
Our next goal is to compute $\tilde{H}^{*,*}(F_{\bullet}\wedge
(BS_l)_+,\zz/l)$ where $S_l$ is the symmetric group and $l$ is a prime
not equal to $char(k)$. 
\begin{lemma}
\llabel{l42} Let $G$ be a finite group and $H$ a subgroup of
$G$. Assume that $[G:H]$ is invertible in the coefficients ring
$R$. Then, for any pointed simplicial sheaf $F_{\BB}$, the map of
motivic cohomology
\begin{equation}
\llabel{splits}
\tilde{H}^{*,*}(F_{\BB}\wedge BG_+,R)\sr \tilde{H}^{*,*}(F_{\BB}\wedge BH_+,R)
\end{equation}
is a split mono and its image is contained in the invariants under the
action of the normalizer of $H$ in $G$.
\end{lemma}
\begin{proof}
We will use the notations established at the beginning of Section
\ref{sec6}. Choose a linear representation $G\sr GL(V)$. The map
(\ref{splits}) is defined by the collection of maps $p_n:\tilde{V}_n/H\sr
\tilde{V_n}/G$ with respect to the identification of Corollary
\ref{limit}. The maps $p_n$ are finite etale of degree $[G:H]$ and the
fundamental cycle on $\tilde{V}_n/H$ over $\tilde{V_n}/G$ defines a
map of freely generated sheaves with transfers
\begin{equation}\llabel{trmap}
p^{\#}:\zz_{tr}(\tilde{V_n}/G)\sr \zz_{tr}(\tilde{V}_n/H)
\end{equation}
The composition of $p^{\#}\zz_{tr}(p)$ is the multiplication by
$[G:H]$. Replacing $F_{\BB}$ by its the standard simplicial resolution
by coproducts of representable sheaves we may assume that terms of
$F_{\BB}$ are coproducts of sheaves of the form $(h_U)_+$. Then,
$$\tilde{H}^{p,q}(F_{\BB}\wedge X_+, A)=Hom_{DM}(N(\zz_{tr}(F_{\BB}))\oo
\zz_{tr}(X),A(q)[p])$$
where $N(-)$ is the normalized chain complex functor from simplicial
sheaves with transfers to complexes of sheaves with transfers. In
particular these groups are functorial in $\zz_{tr}(-)$ which implies
the first claim of the proposition.

The second claim follows from the fact that the normalizer
of $H$ in $G$ acts on $\tilde{V}_n/H$ over $\tilde{V_n}/G$ and the map
$p^{\#}$, being defined by the fundamental cycle, is invariant under
this action.
\end{proof}
Let $\xi_l$ be the vector bundle on
$BS_l$ corresponding to the tautological permutational
representation. Then we have a monomorphism ${\cal O}\sr \xi_l$ and
the quotient $\xi_l/{\cal O}$ is again a vector bundle. Let 
$$d:=e(\xi_l/{\cal O})$$
Assume that $l$ is a prime different from $char(k)$ and that there
exists a primitive l-th root of unity $\zeta$ in $k$. The choice of
$\zeta$ defines a weak equivalence $B{\bf \mu}_l\sr B\zz/l$ and the
inclusion $\zz/l\sr S_l$ gives a map
$$p_{\zeta}:B{\bf \mu}_l\sr B\zz/l\sr BS_l$$
\begin{lemma}
\llabel{restrd} One has:
$$p_{\zeta}^*(d)=-u^{l-1}$$
\end{lemma}
\begin{proof}
The element $p_{\zeta}^*(d)$ is the Euler class of $\xi/{\cal O}$
where $\xi$ corresponds to the regular representation of $\zz/l$ under
our isomorphism $\zz/l\sr \mu_l$. Therefore, we have
$\xi=\oplus_{i=0}^{l-1} L^i$ where $L$ is the line bundle
corresponding to the tautological 1-dimensional representation of
$\mu_l$. By Lemma \ref{l29} and Lemma \ref{l54} we get 
$$p_{\zeta}^*(d)=\prod_{i=1}^{l-1} (iu)=u^{l-1}$$
\end{proof}
\begin{theorem}
\llabel{l10} Let $l$ be a prime and $k$ a field of characteristic not
equal to $l$. There exists a unique class $c\in
H^{2l-3,l-1}(BS_l,\zz/l)$ such that $\delta(c)=d$ and the restriction
of $c$ to $*$ is zero.
\end{theorem}
\begin{proof}
For $l=2$ we have $S_2=\zz/2=\mu_2$ and our result follows from Lemma
\ref{defu}. Assume that $l>2$. The transfer argument shows that to
prove the theorem for $k$ it is sufficient to prove it for a separable
extension of $k$ of degree prime to $l$. In particular, we may assume
that $k$ contains a primitive $l$-th root of unity $\zeta$. 

To prove the existence of $c$ we need to show that $d$ is an
$l$-torsion element in $H^{*,*}(-,\zz)$. For any ring of coefficients
where $(l-1)!$ is invertible, the map $p_{\zeta}$ defines, by Lemma
\ref{l42}, a split injection:
\begin{equation}
\llabel{inj}
H^{*,*}(BS_l)\sr H^{*,*}(B{\bf \mu}_l)^{Aut({\bf \mu}_l)}
\end{equation}
In particular, since (\ref{inj}) is an injection for coefficients in
$\zz$ localized at $l$, it is sufficient to prove that the image of
$d$ in $H^{*,*}(B{\bf \mu}_l,\zz)$ is an l-torsion element. This
follows from Lemma \ref{restrd}.

To show that $c$ is unique, it is sufficient to check that the map
$$\delta:H^{2l-3,l-1}(BS_l,\zz/l)\sr
H^{2l-2,l-1}(BS_l,\zz)$$
is injective. Injectivity of (\ref{inj}) for $\zz/l$-coefficients
implies that it is sufficient to show that the map
\begin{equation}
\llabel{inj2}
\delta:H^{2l-3,l-1}(B{\bf \mu}_l,\zz/l)^{Aut({\bf \mu}_l)}\sr
H^{2l-2,l-1}(B{\bf \mu}_l,\zz)^{Aut({\bf \mu}_l)}
\end{equation}
is injective. Lemma \ref{l52} implies that for $\zz/l$-coefficients
the right hand side of (\ref{inj}) is of the form
\begin{equation}
\llabel{expl}
H^{*,*}(B{\bf \mu}_l,\zz/l)^{Aut({\bf \mu}_l)}=H^{*,*}[[x,y]]/(x^2=0)
\end{equation}
where $x=vu^{l-2}$ and $y=u^{l-1}$. This descriptions shows that
$$H^{2l-3,l-1}(BS_l,\zz/l)^{Aut({\bf \mu}_l)}=\zz/l$$
generated by $x$. We have $\delta(x)=u^{l-1}\ne 0$. Therefore,
(\ref{inj2}) is injective.
\end{proof}
\begin{lemma}
\llabel{restrc}
Let $l$ be a prime such that $char(k)\ne l$ and $\zeta$ be an l-th
root of unity in $k$. Let further $p_{\zeta}:B\mu_l\sr BS_l$ be the
morphism defined by $\zeta$. Then one has:
$$p_{\zeta}^*(c)=-vu^{l-2}$$
\end{lemma}
\begin{proof}
We may assume that $l>2$. Then, the description (\ref{expl}) implies
that $p_{\zeta}^*(c)=avu^{l-2}$ for an element $a\in (\zz/l)^*$. Since
$$\delta(p_{\zeta}^*(c))=p_{\zeta}(\delta(c))=p_{\zeta}(d)=-u^{l-1}$$
and
$$\delta(vu^{l-2})=u^{l-1}$$
we conclude that $a=-1$.
\end{proof}
\begin{theorem}
\llabel{l10new}
For any pointed simplicial sheaf $F_{\bullet}$ over $k$ one has:
\begin{equation}
\llabel{eq10}
\tilde{H}^{*,*}(F_{\bullet}\wedge
(BS_l)_+,\zz/l)=
\left\{
\begin{array}{ll}
\tilde{H}^{*,*}(F_{\bullet},\zz/l)[[c,d]]/(c^2=\tau d +\rho c) &\mbox{\rm for
$l=2$}\\
\tilde{H}^{*,*}(F_{\bullet},\zz/l)[[c,d]]/(c^2=0)&\mbox{\rm for
$l\ne 2$}
\end{array}
\right.
\end{equation}
where:
\begin{enumerate}
\item $\rho$ is the class of $-1$ in $H^{1,1}(k)$
\item $\tau$ is the generator of $H^{0,1}(k,\zz/2)={\bf \mu}_2(k)$.
\end{enumerate}
\end{theorem}
\begin{proof}
For $l=2$ we have $S_2=\zz/2=\mu_2$ and our result follows from
Theorem \ref{l47}. Assume that $l>2$. We need to show that the map
$$\tilde{H}^{*,*}(F_{\BB})[[c,d]]/(c^2=0)\sr
\tilde{H}^{*,*}(F_{\BB}\wedge(BS_l)_+,\zz/l)$$
is an isomorphism. By the transfer argument we may assume that $k$
contains a primitive l-th root of unity. By Lemma \ref{splits} the
homomorphism
$$\tilde{H}^{*,*}(F_{\BB}\wedge (BS_l)_+,\zz/l)\sr
\tilde{H}^{*,*}(F_{\BB}\wedge 
(B\mu_l)_+,\zz/l)^{Aut(\mu_l)}$$
defined by a choice of a primitive root $\zeta$ is a mono. Therefore,
it is sufficient to show that the composition
\begin{equation}
\llabel{inj3} \tilde{H}^{*,*}(F_{\BB})[[c,d]]/(c^2=0)\sr
\tilde{H}^{*,*}(F_{\BB}\wedge (BS_l)_+)\sr
\tilde{H}^{*,*}(F_{\BB}\wedge (B{\bf \mu}_l)_+)^{Aut({\bf \mu}_l)}
\end{equation}
is an isomorphism. The fact that (\ref{inj3}) is an isomorphism
follows from the analog of formula (\ref{expl}) for
$\tilde{H}^{*,*}(F_{\BB}\wedge (B{\bf \mu}_l)_+)^{Aut({\bf \mu}_l)}$
and Lemmas \ref{restrd}, \ref{restrc}.
\end{proof}
\begin{lemma}
\llabel{l1223}
Let $M$ be a line bundle and $l$ a prime not equal to $char(k)$. Then
$P_l(e(M))=e(M)^l+e(M)d$.
\end{lemma}
\begin{proof}
By the transfer argument we may assume that $k$ contains a primitive
l-th root of unity. In view of Lemma \ref{l42} it is sufficient to
prove our equality in the motivic cohomology of $X\times B{\bf
\mu}_l$. By Lemma \ref{l8} we have $P(e(M))=e(M\oo\xi_l)$.  The vector
bundle $\xi_l$ restricted to $B{\bf \mu}_l$ splits into the sum of
line bundles $L^0\oplus L^1\oplus\dots\oplus L^{l-1}$ where $L$ is the
line bundle corresponding to the tautological 1-dimensional
representation of ${\bf \mu}_l$. By Lemmas \ref{l29} and \ref{l54} we
get:
$$e(M\oo\xi_l)=\prod_{i=0}^{l-1} e(M\oo L^i)=\prod_i
(e(M)+ie(L))=e(M)(e(M)^{l-1}-e(L)^{l-1})$$
Since the restriction of $d$ to $B{\bf \mu}_l$ is $-e(L)^{l-1}$ (by
Lemma \ref{restrd}) this finishes the proof.
\end{proof}

\subsection{Symmetry theorem}
Let $G_1$, $G_2$ be two finite groups acting freely on $U_1$ and $U_2$
respectively. Let further $r_i:G_i\sr S_{n_i}$ be permutational
representations of $G_i$, $i=1,2$, $\xi_i$ the corresponding vector
bundles on $U_i/G_i$ of dimension $n_i$ and $L_i\oplus\xi_i\sr {\cal
O}^{N_i}$ inverses of $\xi_i$. Consider the action of $G_1\times G_2$
on $U_1\times U_2$ and let $\xi_1\oo\xi_2$ be the vector bundle on
$(U_1\times U_2)/(G_1\times G_2)$ corresponding to $r_1\times
r_2$. Consider the vector bundle 
$$L_{12}=(L_1\oo\xi_2)\oplus L_2^{N_1}$$
Then
$$L_{12}\oplus(\xi_1\oo\xi_2)=((L_1\oplus \xi_1)\oo\xi_2)\oplus
L_2^{N_1}=\xi_2^{N_1}\oplus L_2^{N_1}={\cal O}^{N_2N_1}$$
i.e. $L_{12}$ is an inverse of $\xi_1\oo\xi_2$. The canonical map
${\cal O}\sr \xi_1$ gives a monomorphism $i:L_1\oplus L_2^{N_1}\sr
L_{12}$. 
\begin{lemma}
\llabel{l5}
The following diagram of pointed sheaves commutes:
\begin{equation}
\llabel{eq5}
\begin{CD}
K_{i}\wedge Th_{L^n_1}\wedge Th_{L_2^{iN_1}} @>\tilde{P}\wedge Id>>
K_{iN_1}\wedge Th_{L_2^{iN_1}}\\ 
@VId\wedge Th(i)VV @VV\tilde{P}V\\
K_i\wedge Th_{L_{12}} @>\tilde{P}>> K_{iN_1N_2}
\end{CD}
\end{equation}
\end{lemma}
\begin{proof}
Direct comparison.
\end{proof}
\begin{proposition}
\llabel{l28}
Let $u$ be a class in $\tilde{H}^{2i,i}(F_{\bullet})$. Then:
$$P_2(P_1(u))e(\xi_1/{\cal O})^{N_1i}=P_{12}(u)e(\xi_1/{\cal
O})^{N_1i}$$
in $\tilde{H}^{2iN_1N_2,iN_1N_2}(F_{\bullet}\wedge (U_1/G_1\times
U_2/G_2)_+)$, where:
$$P_i=P_{G_i,r_i,U_i,L_i}$$
$$P_{12}=P_{G_1\times G_2, r_1\times r_2, U_1\times U_2,L_{12}}$$
\end{proposition}
\begin{proof}
By Lemma \ref{l5} we have
$$\tilde{P}_2(\tilde{P}_1(u))=Th(i)^*\tilde{P}_{12}(u)$$
or, equivalently,
$${P}_2({P}_1(u) t_{L_1^i})
t_{L_2^{iN_1}}=Th(i)^*({P}_{12}(u) t_{L_{12}})$$
By Lemma \ref{l14} we rewrite it as
$${P}_2({P}_1(u) t_{L_1^i})
t_{L_2^{iN_1}}={P}_{12}(u) e((L_1\oo (\xi_2/{\cal O}))^i)
t_{L_1^i} t_{L_2^{iN_1}}$$
By Lemma \ref{l22} we get
$${P}_2{P}_1(u) P_2(t_{L_1^i})
t_{L_2^{iN_1}}={P}_{12}(u) e((L_1\oo (\xi_2/{\cal O}))^i)
t_{L_1^i} t_{L_2^{iN_1}}$$
By Lemmas \ref{l8} and \ref{l14} we get
$${P}_2{P}_1(u) e((L_1\oo (\xi_2/{\cal O}))^i) 
t_{L_1^i} t_{L_2^{iN_1}}={P}_{12}(u) e((L_1\oo
(\xi_2/{\cal O}))^i) t_{L_1^i} t_{L_2^{iN_1}}$$
By Thom isomorphism \ref{l7}, we get
$${P}_2{P}_1(u) e((L_1\oo (\xi_2/{\cal O}))^i)={P}_{12}(u)
e((L_1\oo (\xi_2/{\cal O}))^i)$$
Multiplying both sides by $e((\xi_1\oo (\xi_2/{\cal O}))^i)$ we get by
Lemma \ref{l29}
$${P}_2{P}_1(u) e(\xi_2/{\cal O})^{iN_1}={P}_{12}(u)
e(\xi_2/{\cal O})^{iN_1}$$
\end{proof}
\begin{cor}
\llabel{neededc}
In the notations of the proposition assume that $r_1=l$ is a prime
different from the characteristic of $k$. Then
$$P_2(P_1(u))=P_{12}(u)$$
\end{cor}
\begin{proof}
Replacing $U_2$ by an affine torsor we may assume that it is
affine. Let $G_2\sr GL(V)$ be the linear representation of $G_2$
corresponding to $r_2$ and $\tilde{V}_m$ the open subset of ${\bf
A}(V^{\oplus m})$ where $G$ acts freely. Then, for some $m$, there
exists a $G_2$-equivariant map $U_2\sr \tilde{V}_m$. By Lemma
\ref{l33} it is sufficient to prove the corollary for
$U_2=\tilde{V}_m$ and $G_2=S_{r_2}$. Proposition \ref{l28}
together with Corollary \ref{limit} shows that we have
$$P_2(P_1(u))e(\xi_1/{\cal O})^{N_1i}=P_{12}(u)e(\xi_1/{\cal
O})^{N_1i}$$
on $F_{\bullet}\wedge (U_1/G_1\times BS_l)$. By Theorem \ref{l10new},
multiplication with $e(\xi_1/{\cal O})$ is injective and we conclude
that $P_2(P_1(u))=P_{12}(u)$.
\end{proof}
\begin{lemma}
\llabel{l31}
Let $U$ be a scheme with a free action of $G$ and let $r:G\sr S_n$ be
a permutational representation of $G$. Consider:
$$P=P_{G\times G, r\times r, U\times U-\Delta(U)}:K_i\wedge ((U\times
U-\Delta(U))/(G\times G))_+\sr K_{in^2}$$
Then $P$ is invariant under the permutation of two copies of $U$.
\end{lemma}
\begin{proof}
The action of $G\times G$ on $U\times U-\Delta(U)$ extends to a free
action of the semidirect product $(G\times G)\propto \zz/2$. The
permutational representation $r\times r$ also extends to a
permutational representation of $(G\times G)\propto \zz/2$. Therefore,
$P$ factors through the map
$$K_i\wedge ((U\times
U-\Delta(U))/(G\times G))_+\sr K_i\wedge ((U\times
U-\Delta(U))/(G\times G)\propto \zz/2)_+$$
which implies that it is symmetric. 
\end{proof}
\begin{lemma}
\llabel{l36}
Let $X$ be a smooth variety and $Z$ a smooth subvariety in $Z$ of
codimension $c$. Then for any pointed simplicial sheaf $F_{\bullet}$
the map
$$\tilde{H}^{2i,i}(F_{\bullet}\wedge X_+)\sr \tilde{H}^{2i,i}(F_{\bullet}\wedge
(X-Z)_+)$$
is an isomorphism for $i<c$.
\end{lemma}
\begin{proof}
Follows from Lemma \ref{often}.
\end{proof}
\begin{lemma}
\llabel{l35}
Under the assumptions of Lemma \ref{l31} the morphism
$$P=P_{G\times G, r\times r, U\times U}:K_i\wedge (U\times
U/G\times G)_+\sr K_{in^2}$$
is invariant under the permutation of two copies of $U$.
\end{lemma}
\begin{proof}
Replacing $U$ by an affine torsor we may assume that it is
affine. Let $G\sr GL(V)$ be a faithful liner representation of $G$ and
$\tilde{V}_m$ the open subset of ${\bf A}(V^{\oplus m})$ where $G$
acts freely. Then, for some $m$, there exists a $G$-equivariant map 
$U\sr \tilde{V}_m$. By Lemma \ref{l33} it is sufficient to show that
$P$ is symmetric for $U=\tilde{V}_m$. This follows from Lemma
\ref{l31} and Lemma \ref{l36} since we may choose $m$ such that the
codimension of $\Delta(\tilde{V}_m)$ is larger than $in^2$.
\end{proof}
Corollary \ref{limit} together with Lemma \ref{l33} implies that there
is a well defined morphism:
$$P_l:K_n\wedge (BS_l)_+\sr K_{nl}$$
\begin{theorem}
\llabel{th4}
The composition 
$$K_i\wedge (BS_l)_+\wedge (BS_l)_+\stackrel{P\wedge
Id}{\sr}K_{il}\wedge (BS_l)_+\stackrel{P}{\sr}K_{il^2}$$
is invariant under the permutation of two copies of $BS_l$. 
\end{theorem}
\begin{proof}
Let $\eta$ be the tautological motivic cohomology class of $K_i$.
We need to show that $P(P(\eta))$ is invariant under the permutation
of two copies of $BS_l$. Corollary \ref{neededc} and Corollary
\ref{limit} imply that $P(P(\eta))=P_{12}(\eta)$ where $P_{12}$ is the
power operation corresponding to $S_l\times S_l\sr S_{l^2}$. We
conclude by Lemma \ref{l35} and, again, Corollary \ref{limit}.
\end{proof}

\subsection{Power operations and the Bockstein homomorphism}
We denote by $\beta$ the Bockstein homomorphism 
$$\tilde{H}^{*,*}(-,\zz/l)\sr
\tilde{H}^{*+1,*}(-,\zz/l)$$
which is defined by the short exact sequence of the coefficients
$$0\sr \zz/l\sr \zz/l^2\sr \zz/l\sr 0$$
It has the same properties as the Bockstein homomorphism in the
ordinary cohomology. In particular, we have $\beta\beta=0$ and
for $u\in \tilde{H}^{p,*}$,
\begin{equation}
\llabel{eq2.6.1}
\beta(u v)=\beta(u) v + (-1)^p u \beta(v)
\end{equation}
The goal of this section is to prove Theorem \ref{andb}. The method we use
follows closely the method used to prove an analogous result in
\cite{SE}. 

Let $U$ and $L$ be as in Construction \ref{c2}. The proper
push-forward of cycles defines a ``transfer'' map
$$tr:\uu{Hom}(Th_{U}(L^i),K_{n,R})\sr \uu{Hom}(Th_{U/G}(L^n),K_{n,R})$$
where $\uu{Hom}(-,-)$ denotes the internal Hom-object in the category
of pointed sheaves. Let $\zz_{(l)}$ be the local ring of $l$ in
$\zz$.  Denote by 
$$\Phi=\uu{Hom}(Th_{U/G}(L^n),K_{nl,\zz_{(l)}})/lIm(tr)$$
the pointed sheaf which corresponds to the quotient of the sheaf of
abelian groups $\uu{Hom}(Th_{U/G}(L^n),K_{nl,\zz_{(l)}})$ by the
subgroup of elements of the form $lx$ where $x$ is in the image of the
transfer map.
\begin{lemma}
\llabel{l2.6.3}
The map of pointed sheaves 
$$K_{n,\zz/l}\sr \uu{Hom}(Th_{U/G}(L^n),K_{nl,\zz/l})$$
adjoint to the power operation $\tilde{P}_l$ lifts to a map of pointed
sheaves
$$K_{n,\zz/l}\sr \Phi$$
\end{lemma}
\begin{proof}
The pointed sheaf $K_{n,\zz/l}$ is a quotient sheaf of the sheaf
$K_{n,\zz_{(l)}}$ and the power operation for integral coefficients defines
a map 
$$K_{n,\zz_{(l)}}\sr \uu{Hom}(Th_{U/G}(L^n),K_{nl,\zz_{(l)}})$$
It is sufficient to show that the composition of this map with the
projection to $\Phi$ factors
through $K_{n,\zz/l}$ i.e. that for two cycles $Z_1, Z_2$ be cycles on
$X\times {\bf A}^n$ with integral coefficients such that $Z_1-Z_2$ is
divisible by $l$ the cycle $\tilde{P}_l(Z_1)-\tilde{P}_l(Z_2)$ is in
$l Im(tr)$. Let $Z_i'$ be the cycle on $(X^l\times {\bf A}^{nl}\times
U)/S_l$ whose pull-back to $X^l\times {\bf A}^{nl}\times U$ is
$(p:X^l\times {\bf A}^{nl}\times U\sr X^l\times {\bf
A}^{nl})^*(Z_i^{\oo l})$. It is sufficient to show that $Z_1'-Z_2'$ is
in $l Im(\pi_*)$ where 
$$\pi:X^l\times {\bf A}^{nl}\times U\sr (X^l\times {\bf A}^{nl}\times
U)/S_l$$
is the projection or, equivalently, that 
$$p^*(Z_1^{\oo l})-p^*(Z_2^{\oo l})=l\pi^*\pi_*(Y)$$
for some $Y$. We have $Z_1-Z_2=lW$ and the left hand side can be
rewritten as $p^*((Z_2+lW)^{\oo l})-p^*(Z_2^{\oo l})$. Since any
$S_l$-invariant cycle with coefficients divisible by $l$ is of the
form $\pi^*\pi_*(-)$ it is sufficient to consider the summands in this
expression with coefficients not divisible by $l^2$. They are of the
form $l(Z_2\oo\dots\oo W\oo\dots\oo Z_2)$. The sum of all such cycles
is, up to multiplication by $(l-1)!$, of the form $\pi^*\pi_*(W\oo
(Z_2)^{l-1})$. 
\end{proof}
\begin{lemma}
\llabel{l2.6.4}
Let $u$ be the tautological class in $\tilde{H}^{2n,n}(K_{n,\zz/l})$. Then
$\beta \tilde{P}_l(u)$ lies in the image of the transfer map
$$\tilde{H}^{2nl,nl}(Th_U(L^n)\wedge K_{n,\zz/l},\zz/l)\sr
\tilde{H}^{2nl,nl}(Th_{U/G}(L^n)\wedge K_{n,\zz/l},\zz/l)$$
\end{lemma}
\begin{proof}
Consider the standard simplicial resolution $G_{\bullet}K_{n,R}$ (see
Section \ref{sec3}). Since a cycle with $\zz/l$-coefficients on a
smooth scheme lifts to an integral cycle, the map
$$G_0K_{n,\zz/l}\sr \uu{Hom}(Th_{U/G}(L^n),K_{n,\zz/l})$$
which represents $\tilde{P}_l(u)$ lifts to a map with values in
$\uu{Hom}(Th_{U/G}(L^n),K_{n,\zz})$. Let $\phi$ be such a
lifting. Then $\partial_0\phi-\partial_1\phi$ lands in $lK_{n,\zz}$
and $(\partial_0\phi-\partial_1\phi)/l$ gives a morphism of complexes
$G_*K_{i,\zz/l}\sr \uu{Hom}(Th_{U/G}(L^n), K_{n,\zz/l})$ which defines
$\beta \tilde{P}_l(u)$. Lemma \ref{l2.6.3} implies immediately that
$(\partial_0\phi-\partial_1\phi)/l$ lifts to a morphism with values in 
the image of the transfer map which implies the statement of the
lemma.
\end{proof}
\begin{lemma}
\llabel{l2.6.5}
For any $F_{\bullet}$, the transfer map in cohomology
$$\tilde{H}^{*,*}(F_{\bullet}\wedge Th_{ES_l}(L^n),\zz/l)\sr
\tilde{H}^{*,*}(F_{\bullet}\wedge Th_{BS_l}(L^n),\zz/l)$$
is zero.
\end{lemma}
\begin{proof}
The composition of the transfer map with the restriction map is the
multiplication with the degree of the covering, in our case
$l!$. Hence it is sufficient to show that the restriction map 
$$\tilde{H}^{*,*}(F_{\bullet}\wedge Th_{BS_l}(L^n),\zz/l)\sr
\tilde{H}^{*,*}(F_{\bullet}\wedge Th_{ES_l}(L^n),\zz/l)$$
is surjective. Since motivic cohomology of $ES_l$ are trivial and
class on the right can be written as $u t$ where $u$ is in
$\tilde{H}^{*,*}(F_{\bullet})$ and $t$ is the Thom class. Any such $u t$
is clearly in the image of the restriction map.
\end{proof}
\begin{theorem}
\llabel{andb}\llabel{p2.6.6}
For any $u\in \tilde{H}^{2d,d}$ one has $\beta P(u)=0$.
\end{theorem}
\begin{proof}
By Lemma \ref{l2.6.4} and Lemma \ref{l2.6.5} we have $\beta(P(u)
t)=\beta\tilde{P}(u)=0$. By (\ref{eq2.6.1}) we get
$$\beta(P(u) t)=\beta(P(u)) t$$
and by the Thom isomorphism theorem we conclude that $\beta(P(u))=0$.
\end{proof}

\subsection{Individual power operations: formulas}
In this section we assume that $l$ is a prime different from the
characteristic of $k$. Let $w$ be a class in
$\tilde{H}^{2d,d}(F_{\bullet},\zz/l)$. By Theorem \ref{l10new}, the class
$P_l(w)$ can be written uniquely as a linear combination of the form:
\begin{equation}
\llabel{eq11}
P_l(u)=\sum_{i\ge 0} C_{i+1,d}(w)cd^i + D_{i,d}(w)d^i
\end{equation}
This defines cohomological operations:
$$C_{i,d}:\tilde{H}^{2d,d}(-,\zz/l)\sr
\tilde{H}^{2d+2(d-i)(l-1)+1,d+(d-i)(l-1)}(-,\zz/l)$$
$$D_{i,d}:\tilde{H}^{2d,d}(-,\zz/l)\sr
\tilde{H}^{2d+2(d-i)(l-1),d+(d-i)(l-1)}(-,\zz/l)$$
Below we use $C_{i}$ instead of $C_{i,d}$ and $D_{i}$ instead of
$D_{i,d}$ when no confusion is possible. Recall that we denote by
$\tau$ the generator of $H^{0,1}(k,\zz/2)$ for $char(k)\ne 2$ and by
$\rho$ the class of $-1$ in $H^{1,1}(k,\zz/2)$.
\begin{lemma}
\llabel{car1}\llabel{l61}
Let $u\in \tilde{H}^{2d,d}(F_{\BB})$, $v\in
\tilde{H}^{2d',d'}(F'_{\BB})$. Then for $l$ odd one has: 
$$D_i(u\wedge v)=\sum_{r=0}^{i} D_r(u)\wedge D_{i-r}(v)$$
$$C_{i+1}(u\wedge v)=\sum_{r=0}^{i} C_{r+1}(u)\wedge
D_{i-r}(v)+D_r(u)\wedge C_{i-r+1}(v)$$
and for $l=2$ one has
$$D_i(u\wedge v)=\sum_{r=0}^{i} D_r(u)\wedge D_{i-r}(v) +
\tau\sum_{r=0,\dots,i-1} C_{r+1}(u)\wedge C_{i-r}(v)$$ 
$$C_{i+1}(u\wedge v)=\sum_{r=0}^{i} C_{r+1}(u)\wedge
D_{i-r}(v)+D_r(u)\wedge C_{i-r+1}(v) + \rho C_{r+1}(u)\wedge C_{i-r+1}(v)
$$
\end{lemma}
\begin{proof}
Follows immediately from Lemma \ref{l22} and Theorem \ref{l10new}.
\end{proof}
\begin{lemma}
\llabel{l62}
Let $t\in \tilde{H}^{2,1}(T,\zz/l)$ be the tautological class. Then one has:
\begin{equation}\llabel{eq62}
\begin{array}{l}
C_{i+1}(u\wedge t)=C_i(u)\wedge t\\\\
D_i(u\wedge t)=D_{i-1}(u)\wedge t
\end{array}
\end{equation}
\end{lemma}
\begin{proof}
By Lemma \ref{l40} we have $P(t)=\delta(T\sr Th(\xi_l/{\cal
O}))(t_{\xi})$. In view of Lemma \ref{l14} and the fact that
$d=e(\xi_l/{\cal O})$ we get $P(t)=t\wedge d$ i.e. $C_{i+1}(t)=0$ for
all $i\ge 0$, $D_1(t)=t$ and $D_i(t)=0$ for $i\ne 1$. Applying Lemma
\ref{l61} we get (\ref{eq62}).
\end{proof}
For $u\in \tilde{H}^{2d,d}$ define;
$$P^i(u)=D_{d-i}(u)$$
$$B^i(u)=C_{d-i}(u)$$
By (\ref{eq62}) we have $P^i(u\wedge t)=P^i(u)\wedge t$ and
$B^i(u\wedge t)=B^i(u)\wedge t$. As shown in the proof of Proposition
\ref{redtoone} we can extend $P^i$ and $B^i$ to operations acting on
motivic cohomology groups $\tilde{H}^{p,q}$ for all $p,q$:
$$P^i:\tilde{H}^{p,q}\sr \tilde{H}^{p+2i(l-1),q+i(l-1)}$$
$$B^i:\tilde{H}^{p,q}\sr \tilde{H}^{p+2i(l-1)+1,q+i(l-1)}$$  
For $l=2$ we denote, following the standard convention,
$$Sq^{2i}=P^i$$
$$Sq^{2i+1}=B^i$$
\begin{theorem}
\llabel{vanth}
For any $i<0$ one has $P^i=B^i=0$. 
\end{theorem}
\begin{proof}
This follows from Proposition \ref{neg}.
\end{proof}
\begin{theorem}
\llabel{zero}
One has $P^0=Id$.
\end{theorem}
\begin{proof}
Proposition \ref{zer} implies that $P^0(u)=au$ where $a$ is a
constant. Lemma \ref{l1223} applied to the canonical line bundle on
${\bf P}^1$ implies that $a=1$.
\end{proof}
\begin{lemma}
\llabel{l2.6.8}
One has $\beta B^i=0$ and $\beta P^i=B^i$.
\end{lemma}
\begin{proof}
Follows immediately from Theorem \ref{p2.6.6}, the fact that
$\beta(c)=d$ and the product formula (\ref{eq2.6.1}) for the Bockstein
homomorphism. 
\end{proof}
\begin{proposition}
\llabel{cart}
For $u,v\in \tilde{H}^{*,*}$ and $l\ne 2$ one has:
$$P^i(u\wedge v)=\sum_{r=0}^{i}P^r(u)\wedge P^{i-r}(v)$$
$$B^i(u\wedge v)=\sum_{r=0}^{i}(B^r(u)\wedge
P^{i-r}(v)+P^r(u)\wedge B^{i-r}(v))$$
For $l=2$ we get:
$$Sq^{2i}(u\wedge v)=\sum_{r=0}^{i}Sq^{2r}(u)\wedge Sq^{2i-2r}(v)
+ \tau \sum_{s=0}^{i-1}Sq^{2s+1}(u)\wedge Sq^{2i-2s-1}(v)$$
$$Sq^{2i+1}(u\wedge v)=\sum_{r=0}^{i}(Sq^{2r+1}(u)\wedge
Sq^{2i-2r}(v)+Sq^{2r}(u)\wedge Sq^{2i-2r-1}(v))+$$ 
$$+ \rho \sum_{s=0}^{i-1}Sq^{2s+1}(u)\wedge Sq^{2i-2s-1}(v)$$ 
\end{proposition}
\begin{proof}
Follows immediately from Lemma \ref{l61} and the vanishing result
\ref{vanth}.
\end{proof}
\begin{lemma}
\llabel{power}
For $u\in \tilde{H}^{2n,n}$ one has $P^n(u)=u^l$.
\end{lemma}
\begin{proof}
Follows from Lemma \ref{power0}.
\end{proof}
\begin{lemma}
\llabel{van2}
For $u\in \tilde{H}^{p,q}$ and $n>p-q$, $n\ge q$ one has $P^n(u)=0$.
\end{lemma}
\begin{proof}
Let $i=n+q-p$ and $j=n-q$. Then $\sigma^i_s\sigma^j_t(u)$ is in
$\tilde{H}^{2n,n}$. By Lemma \ref{power} we get 
$$\sigma^i_s\sigma^j_tP^n(u)=(\sigma^i_s\sigma^j_t(u))^l$$
By our assumption $i>0$ and the right hand side is zero because the
diagonal map $S^1_s\sr S^1_s\wedge S^1_s$ is zero in $H_{\bullet}$.
\end{proof}
We will also use the total power operation:
$$R:\tilde{H}^{*,*}\sr \tilde{H}^{*,*}[[c,d,d^{-1}]]/(c^2=\tau d+\rho c)$$
where $\tau=\rho=0$ for $l\ne 2$ and for $l=2$, $\tau$ and $\rho$ are
as in Theorem \ref{l10new}. We define $R$ by the formula
$$R(u)=\sum_{i}(B^{i-1}(u)cd^{-i}+P^i(u)d^{-i})$$
For $l=2$ this becomes
$$R(u)=\sum_i(Sq^{2i-1}(u)cd^{-i}+Sq^{2i}(u)d^{-i})$$
For $u\in \tilde{H}^{2n,n}$ we have $d^{2n}R(u)=P(u)$. Together with Lemma
\ref{l22} this implies that for any $u$ and $v$ one has
$$R(u v)=R(u) R(v)$$
where the right hand side is to be computed in the ring 
$H^{*,*}[[c,d,d^{-1}]]/(c^2=\tau d+\rho c)$.

\subsection{Adem relations}
\begin{lemma}
\llabel{l2.7.1}
Let $d_1, c_1$ be generators of $H^{*,*}(BS_l,\zz/l)$ and $c_2, d_2$
the generators of the ring appearing in the definition of $R(u)$. Then
one has:
\begin{equation}
\llabel{eqf}
R(d_1)=d_1(1-d_1/d_2)^{l-1}
\end{equation}
\begin{equation}
\llabel{eqs}
R(c_1)=(c_1+(d_1/d_2)c_2)(1-d_1/d_2)^{l-2}
\end{equation}
\end{lemma}
\begin{proof}
By transfer argument we may assume that $k$ contains a primitive l-th
root of unity $\zeta$. Consider the map 
$$\phi=p_{\zeta}^*:H^{*,*}(BS_l)\sr
H^{*,*}(B{\bf \mu}_l)$$
defined by $\zeta$. Then $\phi$ is a mono (by Lemma \ref{l42}) and by
Lemmas \ref{restrd},\ref{restrc} one has
$$\phi(d_1)=-u^{l-1}$$
$$\phi(c_1)=-vu^{l-2}$$
where $u$ and $v$ are the generators from Theorem  \ref{l47}. We get:
$$R(u^{l-1})=R(u)^{l-1}=d_2^{1-l}P(u)^{l-1}=d_2^{1-l}(u^l+ud_2)^{l-1}$$
where the last equality holds by Lemma \ref{l1223}. The right hand
side equals to $-\phi(d_1)(1-\phi(d_1)/d_2)^{l-1}$ which implies
(\ref{eqf}). For $R(vu^{l-2})$ we get
$$R(vu^{l-2})=R(v)R(u)^{l-2}=d_2^{2-l}R(v)P(u)=$$
$$=d_2^{2-l}(P^1(v)+B^0(v)c_2+P^0(v)d_2)(u^{l}+ud_2)^{l-2}$$
By Theorem \ref{zero} $P^0(v)=v$, by Lemma \ref{l2.6.8} and since
$\beta(v)=u$, $B^0(v)=u$. By Lemma \ref{van2}, $P^1(v)=0$ and
therefore our expression equals 
$$d_2^{1-l}(u^{l-1}c_2+vu^{l-2}d_2)(u^{l-2}+d_2)^{l-2}=d_2^{1-l}(-\phi(d_1)c_2-\phi(c_1)d_2)(-\phi(d_1)+d_2)^{l-2}$$
This implies (\ref{eqs}).
\end{proof}
\begin{theorem}
\llabel{ademeven}
Let $l=2$ and $0<a<2b$. Then for $a+b=0(mod 2)$ one has:
$$
Sq^aSq^b=\left\{
\begin{array}{ll}
\sum_{j=0}^{[a/2]}
\left(\begin{array}{c}b-1-j\\a-2j\end{array}\right)Sq^{a+b-j}Sq^j&\mbox{\rm
for $a,b$ odd}\\\\
\sum_{j=0}^{[a/2]}\tau^{j mod 2}
\left(\begin{array}{c}b-1-j\\a-2j\end{array}\right)Sq^{a+b-j}Sq^j&\mbox{\rm
for $a,b$ even}
\end{array}
\right.
$$
and for $a+b=1(mod2)$ one has:
$$Sq^aSq^b=\sum_{j=0}^{[a/2]}
\left(\begin{array}{c}b-1-j\\a-2j\end{array}\right)Sq^{a+b-j}Sq^j+\rho^{(j+1)mod
2}S(a,b)$$
where:
$$
S(a,b)=\left\{
\begin{array}{ll}
\left(\begin{array}{c}b-1-j\\a-2j\end{array}\right)Sq^{a+b+j}Sq^{j-1}
&\mbox{\rm for $a$ even, $b$ odd}\\\\
\left(\begin{array}{c}b-1-j\\a-1-2j\end{array}\right)Sq^{a+b+j-1}Sq^{j}
&\mbox{\rm for $a$ odd, $b$ even}
\end{array}
\right.
$$
\end{theorem}
\begin{proof}
Consider the class $P(P(u))$ for $u\in \tilde{H}^{2n,n}$. Denote by $d_1, c_1$
the generators of the cohomology of $BS_l$ appearing when the first
$P$ is applied and by $d_2,c_2$ the generators of the cohomology of
$BS_l$ appearing when the second $P$ is applied. According to the
symmetry theorem \ref{th4} the resulting expression is symmetric with
respect to the exchange of $(d_1,c_1)$ and $(d_2,c_2)$.  We have (to
simplify the notations we sometimes omit $u$ from our expressions):
$$P(u)=\sum_i Sq^{2n-2i-1}cd^i+Sq^{2n-2i}d^i$$
$$P(P(u))=d_2^{2n}R(P(u))=d_2^{2n}\sum_i
(R(Sq^{2n-2i-1})R(c_1)+R(Sq^{2n-2i}))R(d_1)^i=$$
$$=\sum_{i,j}
d_1^i(d_1+d_2)^id_2^{2n-j-i}((Sq^{2j-1}Sq^{2n-2i-1}c_2+ 
Sq^{2j}Sq^{2n-2i-1})(c_1+(d_1/d_2)c_2)+$$
$$+Sq^{2j-1}Sq^{2n-2i}c_2+ 
Sq^{2j}Sq^{2n-2i})$$
Consider the coefficients in this expression at $1$, $c_1$, $c_2$ and
$c_1c_2$. At $1$ we have:
$$\sum_{i,j}d_1^i(d_1+d_2)^id_2^{2n-j-i}(Sq^{2j}Sq^{2n-2i}+\tau d_1
Sq^{2j-1}Sq^{2n-2i-1})$$
At $c_1c_2$ we have:
$$\sum_{i,j}d_1^i(d_1+d_2)^id_2^{2n-j-i}Sq^{2j-1}Sq^{2n-2i-1}$$
At $c_1$ we have:
$$\sum_{i,j}d_1^i(d_1+d_2)^id_2^{2n-j-i}Sq^{2j}Sq^{2n-2i-1}$$
At $c_2$ we have:
$$\sum_{i,j}d_1^{i}(d_1+d_2)^id_2^{2n-j-i}Sq^{2j-1}Sq^{2n-2i}+$$
$$+\sum_{i,j}d_1^{i+1}(d_1+d_2)^id_2^{2n-j-i-1}(\rho Sq^{2j-1}Sq^{2n-2i-1}+
Sq^{2j}Sq^{2n-2i-1})$$
Set $p=i+r$, $q=2n-j-r$. Then coefficient at $d_1^pd_2^q$ is 
$$\sum_i
\left(\begin{array}{c}i\\p-i\end{array}\right)Sq^{4n-2p-2q+2i}Sq^{2n-2i}+ 
\tau\left(\begin{array}{c}i-1\\p-i\end{array}\right)Sq^{4n-2p-2q+2i-1}Sq^{2n-2i+1}
$$
Coefficient at $c_1c_2d_1^pd_2^q$ is
$$\sum_i
\left(\begin{array}{c}i\\p-i\end{array}\right)Sq^{4n-2p-2q+2i-1}Sq^{2n-2i-1}$$
Coefficient at $c_1d_1^pd_2^q$ is
$$\sum_i
\left(\begin{array}{c}i\\p-i\end{array}\right)Sq^{4n-2p-2q+2i}Sq^{2n-2i-1}$$
Coefficient at $c_2d_1^pd_2^q$ is
$$\sum_i
\left(\begin{array}{c}i\\p-i\end{array}\right)Sq^{4n-2p-2q+2i-1}Sq^{2n-2i}+
\left(\begin{array}{c}i\\p-i-1\end{array}\right)Sq^{4n-2p-2q+2i}Sq^{2n-2i-1}$$
$$+ \rho\sum_i\left(\begin{array}{c}i\\p-i-1\end{array}\right)
 Sq^{4n-2p-2q+2i-1}Sq^{2n-2i-1}$$
Consider the coefficient at $c_1c_2d_1^pd_2^q$ where $p=2^s-1$ for
sufficiently large $s$ and $q=x$. For $p$ of this form, the
coefficient $\left(\begin{array}{c}i\\p-i\end{array}\right)$ is
non-zero if and only if $i=p$ (follows from \cite[]{SE}) and we
conclude that our coefficient is $Sq^{4n-2x-1}Sq^{2n-2^{s+1}+1}$. By
symmetry it equals to the coefficient at $c_1c_2d_1^qd_2^p$. Setting
$a=4n-2x-1$, $b=2n-2^{s+1}+1$ and $j=2n-2i-1$ and using the fact that
\begin{equation}
\llabel{us1}
\left(\begin{array}{c}u\\v\end{array}\right)=\left(\begin{array}{c}2u\\2v\end{array}\right)mod\,\,\,2
\end{equation}
we can write
the later as
$$\sum_{j=1\,mod\,2}
\left(\begin{array}{c}2n-1-j\\a-2j-1\end{array}\right)Sq^{a+b-j}Sq^{j}$$
From the standard relation 
\begin{equation}\llabel{us2}
\left(\begin{array}{c}u\\v-1\end{array}\right)=\left(\begin{array}{c}u+1\\v\end{array}\right)+\left(\begin{array}{c}u\\v\end{array}\right)
\end{equation}
and the fact that
$$\left(\begin{array}{c}u\\v\end{array}\right)=0\,mod\,2$$
if $u$ is even and $v$ is odd we get the first of the equalities
stated in the theorem. A very similar argument starting with the
equality between the coefficients at $c_1d_1^pd_2^q$ and
$c_2d_1^qd_2^p$ gives the third equality - the case of even $a$ and
odd $b$. To prove the case when both $a$ and $b$ are even consider the
coefficient at $d_1^pd_2^q$. Consider the second part of this
coefficient i.e. the sum
$$\sum_i\left(\begin{array}{c}i-1\\p-i\end{array}\right)Sq^{4n-2p-2q+2i-1}Sq^{2n-2i+1}=$$
$$=\sum_{j=i-1}\left(\begin{array}{c}j\\p-1-j\end{array}\right)Sq^{4n-2(p-1)-2q+2j-1}Sq^{2n-2j-1}
$$
This is the coefficient at $c_1c_2d_1^{p-1}d_2^q$ which is equal to the
coefficient at $c_1c_2d_1^qd_2^{p-1}$ i.e. to 
$$\sum_{j=i-1}\left(\begin{array}{c}j\\q-j\end{array}\right)Sq^{4n-2(p-1)-2q+2j-1}Sq^{2n-2j-1}=$$
$$=\sum_i\left(\begin{array}{c}i-1\\q-i+1\end{array}\right)Sq^{4n-2p-2q+2i-1}Sq^{2n-2i+1}
$$
Equating our new expression for the coefficient at $d_1^pd_2^q$ with
the old expression for the coefficient at $d_1^qd_2^p$ we get
$$\sum_i\left(\begin{array}{c}i\\p-i\end{array}\right)Sq^{4n-2p-2q+2i}Sq^{2n-2i}=$$
$$=\sum_i\left(\begin{array}{c}i\\q-i\end{array}\right)Sq^{4n-2p-2q+2i}Sq^{2n-2i}+$$ 
$$+\tau\sum_i\left(\left(\begin{array}{c}i-1\\q-i+1\end{array}\right)+
\left(\begin{array}{c}i-1\\q-i\end{array}\right)\right)
Sq^{4n-2p-2q+2i-1}Sq^{2n-2i+1}
$$
Setting again $p=2^s-1$, $a=2n-2q$ and $b=2n-2p$ and using the
standard relations between the binomial coefficients one recovers the
identity for $Sq^aSq^b$ when both $a$ and $b$ are even. Finally, to
get the identity in the case when $a$ is odd and $b$ is even one uses
Lemma \ref{l2.6.8}, the identity for $a$ and $b$ even and the fact
that $\beta(\tau)=\rho$.
\end{proof}
The proof of the following theorem which provides Adem relations for
odd $l$ follows the same line of arguments as the proof of the
corresponding topological fact given in \cite[]{SE}. 
\begin{theorem}
\llabel{ademodd}
For $l>2$ and $0<a<lb$ one has:
$$P^aP^b=\sum_{t=0}^{[a/l]}(-1)^{a+t}\left(\begin{array}{c}(l-1)(b-t)-1\\a-lt\end{array}\right)P^{a+b-t}P^t$$
And for $0\ge a\ge lb$ one has:
$$P^aB^b=\sum_{t=0}^{[a/l]}(-1)^{a+t}\left(\begin{array}{c}(l-1)(b-t)\\a-lt\end{array}\right)B^{a+b-t}P^t+$$
$$+\sum_{t=0}^{[(a-1)/l]}(-1)^{a+t-1}\left(\begin{array}{c}(l-1)(b-t)-1\\a-lt-1\end{array}\right)P^{a+b-t}B^t+$$
\end{theorem}
\begin{proof}
These relations are exactly the same as the Adem relations in the
topological Steenrod algebra for odd coefficients. The proof of these
relations given in \cite[Theorem VIII.1.6]{SE} works in exactly the
same way in the motivic context as in the topological one if one
replaces the reference to \cite[Corollary VIII.1.2]{SE} with the
reference to our symmetry theorem \ref{th4}. The apparent difference
in signs between our situation and the situation of \cite{SE} is
explained by the fact that the image of $d$ in the cohomology of
$B\zz/l$ is $-u^{l-1}=-w_{2l-2}$.
\end{proof}

\subsection{Motivic Steenrod algebra}
Define the motivic Steenrod algebra $A^{*,*}(k,\zz/l)$ as the
subalgebra in the algebra of bistable cohomological operations in the
motivic cohomology with $\zz/l$ coefficients over $k$ generated by
operations $P^i$, $B^i$, $i\ge 0$ and operations of the form $u\mapsto
au$ where $a\in H^{*,*}(k,\zz/l)$. 

Let $I=(\epsilon_0,s_1,\epsilon_1,s_2,\dots,s_k,\epsilon_k)$ be a
sequence where $\epsilon_i\in \{0,1\}$ and $s_i$ are non-negative
integers. Denote by $P^I$ the product
$$P^I=\beta^{\epsilon_0}P^{s_1}\dots P^{s_k}\beta^{\epsilon_k}$$
A sequence $I$ is called admissible if $s_i\ge
ls_{i+1}+\epsilon_i$. Monomials $P^I$ corresponding to admissible
sequences are called admissible monomials. 
\begin{lemma}
\llabel{gen}
Admissible monomials generate $A^{*,*}(k,\zz/l)$ as a left
$H^{*,*}$-module.
\end{lemma}
\begin{proof}
This follows from the Adem relations and the Cartan formula \ref{cart}.
\end{proof}
Our next goal is to show that the admissible monomials are linearly
independent with respect to the left $H^{*,*}$-module structure on
$A^{*,*}$ and, therefore, form a basis of this module. Consider the
submodule $H^{*,>0}A^{*,*}$ in $A^{*,*}$. The Cartan formulas
\ref{cart} imply that its is a two-sided ideal in $A^{*,*}$. Set
$$A^{*,*}_{rig}= A^{*,*}/H^{*,>0}A^{*,*}$$
Using again the Cartan formula one observes that the action of
$A^{*,*}$ on $H^{*,*}(X)$ defines an action of $A^{*,*}_{rig}$ on
$H^{*,*}(X)/H^{*,>0}H^{*,*}(X)$. Theorem \ref{l47} immediately implies
the following result.
\begin{lemma}
\llabel{bmul}
For any $l$ and any $k$ one has:
$$H^{*,*}(B{\bf \mu}_l)/H^{*,>0}H^{*,*}(B{\bf
\mu}_l)=\zz/l[u,v]/(v^2=0)$$
\end{lemma}
Lemma \ref{l2.7.1} implies the following.
\begin{lemma}
\llabel{our2.2}
Let $u$ and $v$ be as in Lemma \ref{bmul}. Then one has the following
equalities in $H^{*,*}(B{\bf \mu}_l)/H^{*,>0}H^{*,*}(B{\bf
\mu}_l)$:
$$
\begin{array}{ll}
\beta(v)=u& P^i(v)=0\,\,\,\, \mbox{\rm for $i>0$}\\
\beta(u^k)=0& P^i(u^k)=
\left(\begin{array}{c}k\\i\end{array}\right)u^{k+i(l-1)}
\end{array}
$$
\end{lemma}
Let $d(I)$ be the degree of an admissible monomial $P^I$ i.e. the
integer such that $P^I$ is an operation from $\tilde{H}^{*,*}$ to
$\tilde{H}^{*+d(I),*+q}$.  
\begin{proposition}
\llabel{inj4} For any $n\ge 0$ there exists $N$ and an element $w$ in
$H^{*,*}((B{\bf \mu}_L)^N)$ such that the elements $P^I(w)$, for all
$I$ such that $d(I)\le n$, are linearly independent with respect to
the left $H^{*,*}$-module structure on $H^{*,*}((B{\bf \mu}_L)^N)$.
\end{proposition}
\begin{proof}
It is sufficient to show that there exists $w$ such that $P^I(w)$ are
linerly independent in $H^{*,*}((B{\bf
\mu}_L)^N)/H^{*,>0}H^{*,*}((B{\bf \mu}_L)^N)$ with respect to the
$\zz/l$-module structure. One starts with Lemma \ref{our2.2} and uses
exactly the same reasoning as in the proof of \cite[Proposition
VI.2.4]{SE}
\end{proof}
The following is an immediate corollary of the proof of Proposition
\ref{linind}.
\begin{cor}
\llabel{linind}
The admissible monomials are linearly
independent with respect to the left $H^{*,*}$-module structure on
$A^{*,*}$.
\end{cor}
Let $l$ be an odd prime. Denote by $A^{*,*}_{top}$ the
$\zz/l$-submodule of $A^{*,*}$ generated by the admissible
monomials. The Adem relations (Theorem \ref{ademodd}) show that
$A^{*,*}_{top}$ is a subring of $A^{*,*}$. Together with Corollary
\ref{linind} they imply that $A^{*,*}_{top}$ is isomorphic to the
usual topological Steenrod algebra with the second grading given by
assigning the weight $(l-1)i$ to $P^i$ and the weight $0$ to $\beta$.
The Cartan formula (Proposition \ref{cart}) shows that the action of
elements of $A^{*,*}_{top}$ on products of motivic cohomology classes
has the same expansion as in topology. Taken together these
observations show that all the standard results about the topological
Steenrod algebra and its dual in the case of odd coefficients
translate without change to the motivic context.  In what follows we
consider both the odd and even coefficients cases but give the proofs
only in the even case where the motivic Steenrod algebra has a 
more complicated structure than its topological counterpart.

Below we denote by $A^{*,*}\oo_{H^{*,*}}A^{*,*}$ the tensor product of left
$H^{*,*}$-modules $A^{*,*}$. The action of $A^{*,*}$ on $\tilde{H}^{*,*}(X)$
is not, in general, $H^{*,*}$-linear. Hence, we can not speak of the
action of $A^{*,*}\oo_{H^{*,*}}A^{*,*}$ on
$\tilde{H}^{*,*}(X)\oo_{H^{*,*}}\tilde{H}^{*,*}(Y)$. However, since the action of
$A^{*,*}$ is $\zz/l$-linear it defines an action of
$A^{*,*}\oo_{\zz/l}A^{*,*}$ on $\tilde{H}^{*,*}(X)\oo_{\zz/l}\tilde{H}^{*,*}(Y)$. If
$x,y$ are two elements of $A^{*,*}\oo_{\zz/l}A^{*,*}$ which become
equal in the tensor product over $H^{*,*}$ then for any $u$ in
$\tilde{H}^{*,*}(X)\oo_{\zz/l}\tilde{H}^{*,*}(Y)$ we have $x(u)=x(v)$ in
$\tilde{H}^{*,*}(X)\oo_{H^{*,*}}\tilde{H}^{*,*}(Y)$. Therefore, for $x$ in
$A^{*,*}\oo_{\zz/l}A^{*,*}$ and $u$ in
$\tilde{H}^{*,*}(X)\oo_{\zz/l}\tilde{H}^{*,*}(Y)$ there is a well defined element
$x(u)$ in $\tilde{H}^{*,*}(X)\oo_{H^{*,*}}\tilde{H}^{*,*}(Y)$.
\begin{lemma}
\llabel{mil1}
For any element $x$ of $A^{*,*}$ there exists a unique element
$$\psi^*(x)=\sum x_i'\oo x_i''$$
of $A^{*,*}\oo_{H^{*,*}}A^{*,*}$ such that for any $X$ and any $u\in
\tilde{H}^{p,*}(X)$, $v\in \tilde{H}^{*,*}(X)$ one has
$$x(u v)=\sum (-1)^{dim(x''_i)p} x'_i(u)
x''_i(v)$$
\end{lemma}
\begin{proof}
Exactly parallel to the proof of \cite[Lemma 1, p.154]{Milnor3} where
Proposition \ref{inj4} is used to prove uniqueness.
\end{proof}
We will need the following lemma below.
\begin{lemma}
\llabel{l2a}
Let $x$ be an element of $A^{*,*}$ and $\psi^*(x)=\sum x_i'\oo
x_i''$. Then for $a\in H^{*,*}$ one has:
$$\psi^*(xa)=\sum x_i'a\oo x_i''=\sum x_i'\oo x_i''a$$
\end{lemma}
\begin{proof}
By uniqueness part of Lemma \ref{mil1} it is enough to check that for
any $u, v\in \tilde{H}^{*,*}(X)$ one has:
$$(\sum x_i'a\oo x_i'')(u\oo v)=(\sum x_i'\oo x_i''a)(u\oo
v)=(xa)(u v)$$
This follows immediately from definitions.
\end{proof}
Since $H^{*,*}$ is not in the center of $A^{*,*}$, the ring structure
on $A^{*,*}\oo_{\zz/l}A^{*,*}$ does not define a ring structure on 
$A^{*,*}\oo_{H^{*,*}}A^{*,*}$. The best we can get in general is an
action of $A^{*,*}\oo_{H^{*,*}}A^{*,*}$ on $A^{*,*}\oo_{\zz/l}A^{*,*}$
with values in $A^{*,*}\oo_{H^{*,*}}A^{*,*}$ given by
$$(u\oo v)(u'\oo v')=uu'\oo vv'$$
We say that an element $f$ of $A^{*,*}\oo_{H^{*,*}}A^{*,*}$ is an
operator-like element if for any two elements $x, y$ of
$A^{*,*}\oo_{\zz/l}A^{*,*}$ which belome equal in the tensor product
over $H^{*,*}$ one has $fx=fy$. For an operator-line element $f$ and
any other element $x$ the product $fx$ is well defined. If $f$ and $g$
are two operator-like elements the product $fg$ is again
operator-line. This shows that operator-like elements form a ring
which we denote by $(A^{*,*}\oo_{H^{*,*}}A^{*,*})_r$.
\begin{lemma}
\llabel{oplike}
For any $x$ in $A^{*,*}$, $\psi^*(x)$ is an operator-like element. The
map $\psi^*:A^{*,*}\sr (A^{*,*}\oo_{H^{*,*}}A^{*,*})_r$ is a ring
homomorphism.
\end{lemma} 
\begin{proof}
Let $y,z$ be two elements of $A^{*,*}\oo_{\zz/l}A^{*,*}$ which become
equal modulo $H^{*,*}$. To check that $\psi^*(x)y=\psi^*(x)z$ it is
sufficient, in view of Proposition \ref{inj4}, to check that for any
$X$ and any $w_1,w_2\in \tilde{H}^{*,*}(X)$ one has
$\psi^*(x)y(w_1\oo w_2)=\psi^*(x)z(w_1\oo w_2)$. Let $c:\tilde{H}^{*,*}(X)\oo
\tilde{H}^{*,*}(X)\sr \tilde{H}^{*,*}(X)$ be the cup product. Then by definition of
$\psi$ we have 
$$\psi^*(x)y(w_1\oo w_2)=x(c(y(w_1\oo w_2)))$$
$$\psi^*(x)z(w_1\oo w_2)=x(c(z(w_1\oo w_2)))$$
Our assumption on $y,z$ implies that $c(y(w_1\oo w_2))=c(z(w_1\oo
w_2))$.

To prove that $\psi^*$ is a ring homomorphism we have to check that
for $x,y\in A^{*,*}$ and $w_1,w_2\in \tilde{H}^{*,*}(X)$ we have 
$$\psi(xy)(w_1\oo w_2)=\psi(x)(\psi(y)(w_1\oo w_2))$$
This follows immediately from definitions.
\end{proof}
\begin{lemma}
\llabel{ascom}
The comultiplication map $\psi^*$ is associative and commutative.
\end{lemma}
\begin{proof}
The associativity follows immediately from the definition. To verify
commutativity it is sufficient, in view of Lemma \ref{oplike}, to
verify it on generators of the algebra $A^{*,*}$ i.e. on operations
$P^i$ and $B^i$. For this operations commutativity follows directly
from the Cartan formulas (Proposition \ref{cart}).
\end{proof}

\subsection{Structure of the dual to the motivic Steenrod algebra}

Let $f:A^{*,*}\sr H^{*,*}$ be a homomorphism of left
$H^{*,*}$-modules.  Such a homomorphism is said to be homogeneous of
bidegree $(p,q)$ if for any $i,j\ge 0$ it takes $A^{i,j}$ to
$H^{i-p,j-q}$. We denote by $A_{*,*}$ the ``bigraded dual'' to
$A^{*,*}$ i.e. the group of the left $H^{*,*}$-module maps from
$A^{*,*}$ to $H^{*,*}$ which are finite sums of homogeneous maps.

Let $P^I$
be the basis of admissible monomials in $A^{*,*}$ and $\theta(I)^*$
the dual basis in $A_{*,*}$. An element $x$ from $A_{p,q}$ can be
written uniquely as a sum of the form
\begin{equation}
\llabel{fsum}
x=\sum a_I \theta(I)^*
\end{equation}
where $\theta(I)^*\in A^{r,s}$ and $a_I\in H^{i,j}$ such that $p=r-i$,
$q=s-j$. Since $H^{i,j}=0$ for $i>j$ we get
$$r-s=p-q + i -j\le p-q$$
Since for any $I$, $r\ge 2s$ and $s\ge 0$ we further have
$$2s\le r\le p-q+s$$
which implies:
$$0\le s\le p-q;\,\,\,\,\, 2s\le r\le 2(p-q)$$
Therefore, the sum (\ref{fsum}) is always finite. This fact implies in
particular that taking the bigraded dual to $A_{*,*}$ we get back the
original $A^{*,*}$. It also implies that $A_{*,*}$ is a free $H^{*,*}$
module. 

We have a homomorphism $H^{*,*}\sr A_{*,*}$ which takes $a\in
H^{p,q}$ to the map $\phi\mapsto a\phi(1)$ which lies in $A_{-p,-q}$ and
which we also denote by $a$. The homomorphism $\psi^*$ defines a
homomorphism $\phi_*:A_{*,*}\oo A_{*,*}\sr A_{*,*}$. Lemma \ref{ascom}
immediately implies the following result.
\begin{proposition}
\llabel{isaring}
The homomorphism $\phi_*$ makes $A_{*,*}$ into an associative ring
which is graded commutative with respect to the first grading.
\end{proposition}
Let $e^i\in A^{p_i,q_i}$ be a basis of $A^{*,*}$ over $H^{*,*}$ such
that:
\begin{enumerate}
\item $e_0=1$ and $q_i>0$ for $i>0$
\item $p_i\ge 2q_i$ 
\item for any $q$ there are only finitely many $i$ with $q_i\le q$.
\end{enumerate}
An example of such a basis is given by the basis of admissible
monomials. Let $e_i$ be the dual basis of $A_{*,*}$. Let $X$ be a
smooth scheme over $k$. Then $H^{p,q}(X)=0$ for $p>q+dim(X)$ and
therefore for any $w\in H^{*,*}(X)$ there is only finitely
many $i$'s such that $e^i(w)\ne 0$. We can define a map
$$\lambda^{*,*}:H^{*,*}(X)\sr H^{*,*}(X)\oo_{H^{*,*}}A_{*,*}$$
by the formula
$$\lambda^{*,*}(w)=\sum e^i(w)\oo e_i$$
The following lemma is straightforward.
\begin{lemma}
\llabel{donot}
The map $\lambda^{*,*}$ is a ring homomorphism which does not depend
on the choice of $e^i$.
\end{lemma}
Note that in the case when $X=Spec(k)$ this homomorphism does not
coincide with the ``scalar'' map $H^{*,*}\sr {\cal A}_{*,*}$ which is
described above since the action of $e_i$'s on $H^{*,*}$ may be
nontrivial. In particular $\lambda^{*,*}$ is not a $H^{*,*}$-module
homomorphism.

If we consider $B{\bf \mu}_l$ as a colimit of smooth schemes we can
write a formal analog of $\lambda^{*,*}$. In particular for the
canonical generators $u$ and $v$ we get:
$$\lambda^{*,*}(v)=\sum_{i=0}^{\infty}(vu^i\oo x_i +
u^i\oo y_i)$$
$$\lambda^{*,*}(u)=\sum_{i=0}^{\infty}(vu^i\oo x'_i +
u^i\oo y'_i)$$ 
where $x_i,y_i,x_i',y_i'$ are some well defined elements in ${\cal
A}_{*,*}(k,\zz/l)$. We denote:
$$\xi_i=y'_{l^{i}}\in {\cal A}_{2(l^i-1),l^i-1}(k,\zz/l)$$
$$\tau_i=y_{l^i}\in {\cal A}_{2l^i-1,l^i-1}(k,\zz/l)$$
Since for any basis of $A^{*,*}$ such that $e^0=1$ we have
$dim(e^i)>0$ for $i\ne 0$ we conclude that $\xi_0=1$. Denote by $M_k$
the monomial $P^{l^{k-1}}P^{l^{k-2}}\dots P^1$. 
\begin{lemma}
\llabel{ourl5}
We have in $H^{*,*}(B{\bf \mu}_l)$:
\begin{equation}
\llabel{l5eq}
M_k(u)=u^{l^k}\,\,\,\, M_k\beta(v)=u^{l^k}
\end{equation}
If $f$ is any monomial in $P^i$ and $\beta$ other than $M_k$ (resp. other
than $M_k\beta$) then $f(u)=0$ (resp. $f(v)=0$).
\end{lemma}
\begin{proof}
The equations (\ref{l5eq}) follow from Lemma \ref{power} and the fact
that $\beta(v)=u$. The other statement follows from Lemma \ref{van2},
Lemma \ref{l1223} and multiplicativity of $P$.
\end{proof}
Taking the basis of admissible monomials to compute $\lambda^{*,*}$
and using Lemma \ref{ourl5} we conclude that 
\begin{equation}
\llabel{eqtoth}
\lambda^{*,*}(v)=v\oo 1+\sum_{i=0}^{\infty}u^{l^i}\oo \tau_i
\end{equation}
$$\lambda^{*,*}(u)=u\oo 1
+\sum_{i=0}^{\infty}u^{l^i}\oo \xi_i.$$ 
For an element $\phi$ in $A^{*,*}$ and an element $\psi$ in $A_{*,*}$
let $\langle\psi,\phi\rangle\in H^{*,*}$ be the value of $\psi$ on
$\phi$. Then we 
have:
\begin{equation}
\llabel{eqtoth2}
\phi(u)=\langle\xi_0,\phi\rangle u+\sum_{i}\langle \tau_i,\phi\rangle v^{l^i}
\end{equation}
\begin{equation}
\llabel{eqtoth2.2}
\phi(v)=\sum_i\langle \xi_i,\phi\rangle v^{l^i}
\end{equation}
Let $I$ be a sequence of the form
$(\epsilon_0,r_1,\epsilon_1,r_2,\dots,)$ where $\epsilon_i\in \{0,1\}$,
$r_i\ge 0$ are nonnegative integers and $I$ has only finitely many
nontrivial terms. Following \cite{Milnor3} we set:
$$\omega(I)=\tau_0^{\epsilon_0}\xi_1^{r_1}\tau_1^{\epsilon_1}\xi_2^{r_2}\dots$$
in  ${\cal A}_{p(I),q(I)}(k,\zz/l)$ where
$$p(I)=\epsilon_0+\sum_{i\ge 1}(\epsilon_i(2l^i-1)+2r_i(l^i-1))$$
$$q(I)=\sum_{i\ge 1}(\epsilon_i+r_i)(l^i-1).$$
and
$$\theta(I)=\beta^{\epsilon_0}P^{s_1}\beta^{\epsilon_1}P^{s_2}\dots$$
where 
$$s_n=\sum_{i\ge n}(\epsilon_i+r_i)l^{i-n}.$$
Simple computation shows that $\theta(I)$ belongs to ${\cal
A}^{p(I),q(I)}(k,\zz/l)$. 

In the following theorem we consider, following \cite{Milnor}, the
lexicographical order on the set of sequences $I$ such that
$(1,2,0,\dots)<(0,0,1,\dots)$.
\begin{theorem}
\llabel{howlong0} Then, for $I<J$ one has $\langle
\theta(I),\omega(J)\rangle=0$ and for $I=J$ one has $\langle
\theta(I),\omega(J)\rangle=\pm 1$
\end{theorem}
\begin{proof}
The value of $\langle \theta(I),\omega(J)\rangle$ is a homogeneous
element of $H^{*,*}$ of degree zero. Hence, it is sufficient to show
that the image of $\langle \theta(I),\omega(J)\rangle$ in
$H^{*,*}/H^{*,>0}$ is $0$ or $\pm 1$ depending on whether $I<J$ or
$I=J$. This is done using the action of $A^{*,*}/A^{*,>0}$ on
$H^{*,*}(B{\bf \mu}_l)/H^{*,>0}H^{*,*}(B{\bf \mu}_l)$ described in
Lemma \ref{our2.2} in exactly the same way as in the proof of
\cite[Lemma 8, p.160]{Milnor3}.
\end{proof}
\begin{cor}
\llabel{howlong} Elements $\omega(I)$ (resp. $\theta(I)$) form a basis
of the $H^{*,*}$-module ${\cal A}_{*,*}(k,\zz/l)$ (resp. ${\cal
A}^{*,*}(k,\zz/l)$).
\end{cor}
\begin{proof}
Elements $\theta(I)$ are exactly the admissible monomials. They form a
basis of $A^{*,*}$ by Lemma \ref{gen} and Corollary
\ref{linind}. The fact that elements $\omega(I)$ form a basis of
$A_{*,*}$ follows now from Theorem \ref{howlong0}.
\end{proof}
\begin{theorem}
\llabel{dual}
The (graded commutative) algebra ${\cal A}^{*,*}(k,\zz/l)$ over $H^{*,*}$
is canonically isomorphic to the (graded commutative) algebra with generators
$$\tau_i\in A_{2l^i-1,l^i-1}$$
$$\xi_i\in A_{2l^i-2,l^i-2}$$
and relations
\begin{enumerate}
\item $\xi_0=1$
\item $\tau_i^2=\left\{
\begin{array}{ll}
0&\mbox{\rm for $l\ne 2$}\\
\tau\xi_{i+1}+\rho\tau_{i+1}+\rho\tau_0\xi_{i+1}&\mbox{\rm for $l=2$}
\end{array}
\right.
$
\end{enumerate}
\end{theorem}
\begin{proof}
We already know by Corollary \ref{howlong} that ${\cal
A}_{*,*}(k,\zz/l)$ has a basis which consist of monomials in
$\xi_i,\tau_i$ which are of degree $\le 1$ in each $\tau_i$. The
relation $\tau_i^2=0$ for odd $l$ is a corollary of graded
commutativity. Thus we may assume that $l=2$ in which case we have
only to show that
$\tau_i^{2}=\tau\xi_{i+1}+\rho\tau_{i+1}+\rho\tau_0\xi_{i+1}$.  The
required relation follows immediately from (\ref{eqtoth}), the
multiplicativity of $\lambda^{*,*}$ and the relation
$v^2=\tau u+\rho v$ in the motivic cohomology ring of
$B\mu_2$.
\end{proof}
\begin{lemma}
\llabel{forq1}
For any $\phi\in A^{*,*}$ one has:
\begin{equation}
\llabel{eqtoth3}
\phi(u^{l^j})=\sum_i\langle \xi_{i}^{l^j},\phi\rangle u^{l^{i+j}}
\end{equation}
\end{lemma}
\begin{proof}
Let $x_1,x_2$ be any elements of $A_{*,*}$. If $\psi^*(\phi)=\sum
\phi_k'\oo\phi_k''$ we have, by definition of 
product in $A_{*,*}$:
$$\langle  x_1x_2,\phi\rangle =\sum \langle  x_1,\phi_k'\rangle\langle  x_2,\phi_k''\rangle$$
This implies by induction starting with (\ref{eqtoth2}) that for any
$n$ one has
$$\phi(u^n)=
\sum_{(j_1,\dots,j_n)}\langle \xi_{j_1}\dots\xi_{j_n},\phi\rangle u^{l^{j_1+\dots+j_n}}$$ 
For $n=l^j$ all the terms except for the ones which show up in the
right hand side of (\ref{eqtoth3}) cancel out since we work with
$\zz/l$-coefficients. 
\end{proof}
\begin{proposition}
\llabel{copprop} Let $\phi,\psi$ be elements of $A^{*,*}$ such that
$$\langle \xi_k,\psi\rangle , \langle \tau_k,\psi\rangle \in
\zz/l\subset H^{*,*}$$
Then one has:
$$\langle \xi_k,\phi\psi\rangle=\sum_i\langle \xi_{k-i}^{l^i},\phi\rangle\langle \xi_i,\psi\rangle$$
$$\langle \tau_k,\phi\psi\rangle=\langle \tau_k,\phi\rangle\langle \xi_0,\psi\rangle+
\sum_i\langle \xi_{k-i}^{l^i},\phi\rangle\langle \tau_i,\psi\rangle$$
\end{proposition}
\begin{proof}
We have by (\ref{eqtoth2})
$$\phi\psi(u)=\sum \langle \xi_i,\phi\psi\rangle u^{l^i}$$
On the other hand by (\ref{eqtoth2}) and (\ref{eqtoth3}) we get:
$$\phi\psi(u)=\phi(\sum
\langle \xi_i,\psi\rangle
u^{l^i})=\sum_{i,j}\langle \xi_{j}^{l^i},\phi\rangle\langle \xi_i,\psi\rangle
u^{l^{i+j}}
$$
Comparing coefficients at powers of $u$ we get the first of the
required equalities. To get the second one we write
$$\phi\psi(v)=\langle \xi_0,\phi\psi\rangle
v+\sum_{i}\langle \tau_i,\phi\psi\rangle u^{l^i}$$
by (\ref{eqtoth2}). And by (\ref{eqtoth2}) and (\ref{eqtoth3}) we get:
$$\phi\psi(v)=\langle \xi_0,\psi\rangle\phi(v)+\sum_{i}\langle \tau_i,\psi\rangle\phi(u^{l^i})=$$
$$=\langle \xi_0,\phi\rangle\langle \xi_0,\psi\rangle
v+\sum_{i}\langle \tau_i,\phi\rangle\langle \xi_0,\psi\rangle u^{l^i}+
\sum_{i,j}\langle \xi_{j}^{l^i},\phi\rangle\langle \tau_i,\psi\rangle
u^{l^{i+j}}$$ 
Comparing coefficients we get the second equality.
\end{proof}
Now we can describe the dual to the ring structure on $A^{*,*}$. We
have two $H_{*,*}$ module structures on $A_{*,*}$. The first one, the
left module structure which we used all the time, is given by
$$a*_l\phi(x)=\phi(ax)=a\phi(x)$$
where $\phi\in A_{*,*}$, $a\in H^{*,*}$ and $x\in A^{*,*}$. The other
one is the right module structure given by
$$\phi*_ra(x)=\phi(xa)$$
Lemma \ref{l2a} implies that $(\phi\phi')*_r a=\phi(\phi'*_r a)$ and,
in particular, that $\phi*_r a=\phi(1*_r a)$. The map $a\mapsto 1*_r
a$ coincides with the map $\lambda^{*,*}$ for $X=Spec(k)$ and we
denote it by $\lambda$. Therefore, the two module structures are
defined by two ring homomorphisms $a\mapsto a\xi_0$ and $a\mapsto
\lambda(a)$ from $H^{*,*}$ to $A_{*,*}$.

Denote by $A_{*,*}\oo_{r,l}A_{*,*}$ the tensor product with the
property 
$$(\phi*_ra)\oo\psi=\phi\oo(a*_l\psi).$$ 
Similarly, denote by $A^{*,*}\oo_{r,l}A^{*,*}$ the tensor product with
the property 
$$xa\oo y=x\oo ay.$$
The following lemma is taken from \cite[Lemma 3.3]{Boardman}.
\begin{lemma}
\llabel{boardman} Let $f,g$ be elments of $A_{*,*}$ and $x,y$ elements
of $A^{*,*}$. The formula:
\begin{equation}
\llabel{beq}
\langle \theta(f\oo g), x\oo y\rangle =(-1)^{deg(g)deg(m)}\langle f,x\langle g,y\rangle \rangle 
\end{equation}
defines an isomorphism 
$$\theta:A_{*,*}\oo_{r,l}A_{*,*}\sr (A^{*,*}\oo_{r,l}A^{*,*})^*$$
where the upper star on the right hand side denotes the bigraded dual
of left $H^{*,*}$-module maps from $A^{*,*}\oo_{H^{*,*}}A^{*,*}$ to
$H^{*,*}$.
\end{lemma}
\begin{proof}
One verifies easily that $\theta$ is indeed well defined by
(\ref{beq}). To prove that it is an isomorphism consider the basis
$\omega(I)$ in $A_{*,*}$ and let $\omega(I)^*$ be the dual basis in
$A^{*,*}$. The elements $\omega(I)\oo \omega(J)$ clearly generate
$A_{*,*}\oo_{r,l}A_{*,*}$ as a left $H^{*,*}$-module. The image of
$\omega(I)\oo \omega(J)$ with respect to $\theta$ is the functional
which equals one on $\omega(I)^*\oo \omega(J)^*$ and zero on all other
elements of the form $\omega(I')^*\oo \omega(J')^*$. This implies that
$\omega(I)\oo \omega(J)$ are linearly independent in
$A_{*,*}\oo_{r,l}A_{*,*}$ and hence form a basis of this left
$H^{*,*}$-module. It also implies that $\theta$ maps this basis to a
basis of $(A^{*,*}\oo_{H^{*,*}}A^{*,*})^*$ and therefore $\theta$ is
an isomorphism.
\end{proof}
Composing the dual to the multiplication map $A^{*,*}\oo A^{*,*}\sr
A^{*,*}$ with $\theta$ we get a map
\begin{equation}
\llabel{copr2}
\psi_*:A_{*,*}\sr A_{*,*}\oo_{r,l}A_{*,*}
\end{equation}
By construction, the map $\psi_*$ is defined by the property that
$\psi_*(f)=\sum f_i'\oo f_i''$ and for any $x,y\in A^{*,*}$ one has:
\begin{equation}
\llabel{copr3}
\sum \langle f_i', x\langle f_i'',y\rangle \rangle =\langle f,xy\rangle 
\end{equation}
\begin{lemma}
\llabel{isaring2}
The map (\ref{copr2}) is a ring homomorphism.
\end{lemma}
\begin{proof}
It follows from a direct computation and Lemma \ref{l2a}.
\end{proof}
In view of Lemma \ref{isaring2} and Theorem \ref{dual}, the map $\psi_*$ is  
completely determined by its values on the generators
$\xi_i,\tau_i$. 
\begin{lemma}
\llabel{ongen}
One has:
$$\psi_*(\xi_k)=\sum_{i=0}^k\xi^{l^i}_{k-i}\oo\xi_i$$
$$\psi_*(\tau_k)=\sum_{i=0}^k\xi^{l^i}_{k-i}\oo\tau_i+\tau_k\oo 1$$
\end{lemma} 
\begin{proof}
Follows from Proposition \ref{copprop} and the formula (\ref{copr3}).
\end{proof}
\begin{remark}\llabel{hopfal}\rm
The rings $H^{*,*}$ and $A_{*,*}$, two homomorphisms $H^{*,*}\sr
A_{*,*}$, the homomorphism $A_{*,*}\sr H^{*,*}$ which takes $\tau_i$
for $i\ge 0$ and $\xi_i$ for $i>0$ to zero and the homomorphism
$\psi_*$ form together a Hopf algebroid ${\cal H}(k,\zz/l)$. We can
not give its complete description because we do not know the structure
of $H^{*,*}$ and the explicit form of the homomorphism $\lambda$ which
involves the action of the reduced power operations and the Bockstein
in $H^{*,*}$. One can easily see however that these are the only two
pieces of information missing. In the case when $l>2$ we have a
coaction of the topological dual Steenrod algebra $A_{*}(l)$ (given the
second grading in the way explained above) on $H^{*,*}$ and ${\cal
H}(k,\zz/l)$ is the twisted product of $A_*(l)$ and $H^{*,*}$ with
respect to this coaction. For $l=2$ consider the Hopf algebroid ${\cal
H}(2)$ over $\zz/2$ defined as follows:
\begin{description}
\item[Ring of objects] is $\zz/2[\rho,\tau]$ where $deg(\rho)=(-1,-1)$
and $deg(\tau)=(-1,0)$
\item[Ring of morphisms] is 
$$\zz/2[\rho,\tau,\tau_i,\xi_{i+1}]_{i\ge
0}/(\tau_i^2-\tau\xi_{i+1}-\rho\tau_{i+1}-\rho\tau_0\xi_{i+1})_{i\ge 0}$$
\item[Coface maps] are given by
$$d_0(\rho)=\rho\,\,\,\,\,d_0(\tau)=\tau$$
$$d_1(\rho)=\rho\,\,\,\,\,d_1(\tau)=\tau+\rho\tau_0$$
\item[Codegeneracy map] is given by
$$s_0(\rho)=\rho\,\,\,\,\, s_0(\tau)=\tau$$
$$s_0(\tau_i)=0\,\,\,\mbox{\rm for $i\ge 0$}$$
$$s_0(\xi_i)=0\,\,\,\mbox{\rm for $i>0$}$$
\item[Co-composition] is given by
$$\psi_*(\rho)=\rho\oo 1=1\oo \rho$$
$$\psi_*(\tau)=\tau\oo 1=1\oo\tau +\rho\tau_0\oo 1$$
$$\psi_*(\xi_k)=\sum_{i=0}^k\xi^{l^i}_{k-i}\oo\xi_i$$
$$\psi_*(\tau_k)=\sum_{i=0}^k\xi^{l^i}_{k-i}\oo\tau_i+\tau_k\oo 1$$
\end{description}
Note that our formulas imply that ${\cal H}(2)$ is, in fact, a Hopf
algebroid over $\zz/2[\rho]$. This Hopf algebroid co-acts on
$H^{*,*}$ and ${\cal H}(k,\zz/2)$ is the twisted product of ${\cal
H}(2)$ and $H^{*,*}$ with respect to this coaction.
\end{remark}

\subsection{Operations $\rho(E,R)$ and their properties}
\llabel{sec13}
Let $R=(r_1,r_2,\dots)$ be a sequence of non-negative integers which
are almost all zero and $E=(\epsilon_0,\epsilon_1,\dots)$ a sequence
of zeros and ones which are almost all zeros. Corollary \ref{howlong}
implies that elements of the form
$$\tau(E)\xi(R):=\prod_{i\ge 0}\tau_i^{\epsilon_i}\prod_{j\ge
1}\xi_{j}^{r_j}$$ 
form a basis of $A_{*,*}$ over $H^{*,*}$. Let $\rho(E,R)$ be the dual
basis of $A^{*,*}$. In particular, we define
$${\cal P}^R=\rho(\uu{0},R)$$
$$Q(E):=\rho(E,\uu{0})$$
and 
$$Q_i=Q(0,\dots,0,1,0,\dots)$$
where $1$ is on the i-th place i.e. $Q_i$ is the dual to $\tau_i$. 

If $l$ is odd operations $\rho(E,R)$ and, in particular, $Q(E)$,
${\cal P}^R$ and $Q_i$, have the same properties as their topological
counterparts defined in \cite{Milnor3}. In what follows we assume that
$l=2$.

\begin{lemma}
\llabel{dt}
${\cal P}^{(n,0,\dots,0,\dots)}=P^n$
\end{lemma}
\begin{proof}
Follows immediately from Theorem \ref{dual}.
\end{proof}
\begin{proposition}
\llabel{com0} $\rho(E,R)=Q(E){\cal P}^R$.
\end{proposition}
\begin{proof}
We have to compute the pairing of $Q(E){\cal P}^R$ with
$\tau(E')\xi(R')$ and show that it is $1$ for
$E'=E$, $R'=R$ and zero otherwise. By
(\ref{copr3}) we have:
\begin{equation}
\llabel{expr1}
\langle \tau(E')\xi(R'),Q(E){\cal P}^R\rangle =\sum
\langle f_i',Q(E)\langle f_i'',{\cal P}^R\rangle \rangle 
\end{equation}
where 
$$\sum f_i'\oo
f_i''=\psi_*(\tau(E')\xi(R'))=\psi_*(\tau(E'))
\psi_*(\xi(R'))$$ 
We can choose our representation $\sum f_i'\oo f_i''$ such that
$f_i''$ are of the form $\rho(E'',R'')$ and, in particular, $\langle
f_i'',{\cal P}^R\rangle $ are in $\zz/2$. Then, the expression (\ref{expr1})
depends only on the class of $\psi_*(\tau(E')\xi(R'))$
in $A_{*,*}/I\oo_{r,l}A^{*,*}$ where $I$ is the ideal generated by
$\xi_i$ for $i>0$. In this quotient ring we have:
$$\psi_*(\xi_k)=1\oo \xi_k$$
$$\psi_*(\tau_k)=1\oo\tau_k +\tau_k\oo 1$$
This easily implies that (\ref{expr1}) is non-zero if and only if
$E=E'$ and $R=R'$.
\end{proof}
Let $\rho$ be, as before, the class of $-1$ in $H^{1,1}$. Denote by
$B^{*,*}$ the $\zz/2[\rho]$-submodule in $A^{*,*}$ generated by
elements of the form $Q(E)$. Let $B_{*,*}$ be the dual of
$B^{*,*}$. Then
$$B_{*,*}\oo_{\zz/2[\rho]}H^{*,*}=A_{*,*}/(\{\xi_i\}).$$
Lemma \ref{ongen} implies that for $f\in I=(\{\xi_i\})$, one has
$$\psi_*(f)\in I\oo A_{*,*}+A_{*,*}\oo I$$
and, therefore, $\psi_*$ defines a comultiplication on
$B_{*,*}\oo_{\zz/2[\rho]}H^{*,*}$ which takes $\tau_k$ to $\tau_k\oo
1+1\oo \tau_k$. From this one easily gets the following result.
\begin{proposition}
\llabel{bstr} As a  $\zz/2[\rho]$-algebra, $B_{*,*}$ is of the form 
$$B_{*,*}=\zz/2[\rho][\tau_0,\dots,\tau_n,\dots]/(\tau_i^2=\rho\tau_{i+1})$$
The map $\psi_*$ defines a Hopf algebra structure on $B_{*,*}$ over
$\zz/2[\rho]$, satisfying:
$$\phi_*(\tau_i)=\tau_i\oo 1+1\oo \tau_i.$$
\end{proposition}
Dualizing we get the following theorem on the structure of $B^{*,*}$.
\begin{proposition}
\llabel{bstrd}
As a $\zz/2[\rho]$-algebra, $B^{*,*}$ is the exterior algebra with
generators $Q_i$. For $E=(\epsilon_0,\dots,\epsilon_n)$ one
has $Q(E)=\prod Q_i^{\epsilon_i}$. The coproduct is given
on $Q_i$'s by
$$\psi^*(Q_i)=1\oo Q_i + Q_i\oo
1+\rho\sum_{E,E'}
c_{E,E'}Q(E)\oo
Q(E')$$ 
where $E,E'$ run through sequences of zeros and ones which are almost
all zeros and $c_{E,E'}$ are elements of $H^{*,*}$.
\end{proposition}
The following three results complete the proof of all the properties of
operations $Q_i$ used in \cite{MC}.
\begin{lemma}
\label{previous}
$Q_0=\beta$.
\end{lemma}
\begin{proof}
Since operations $\rho(E,R)$ form a basis we can write $\beta$ as a
formal linear combination $\sum a_{E,R}\rho(E,R)$. Since the weight of
$\beta$ is zero this implies that $\beta=cQ_0$ for $c\in \zz/l$. Since
$\beta(u)=v$, formula (\ref{eqtoth2}) implies that $c=1$.
\end{proof}
\begin{proposition}
\llabel{com1} Let $i\ge 1$ and $q_i={\cal P}^{0,\dots,0,1,0,\dots}$ be
the dual to $\xi_i$. Then one has
$$Q_{i}=[Q_0,q_i]$$
\end{proposition}
\begin{proof}
We have to show that $q_iQ_0=Q_0q_i+Q_{i}$ i.e. that the only
monomials which pair non-trivially with $q_iQ_0$ are $\tau_0\xi_i$ and
$\xi_{i}$ and that for those monomials the pairing gives $1$. Using
formula (\ref{copr3}) we see that it is sufficient to show that the
only monomials $M=\tau(E)\xi(R)$ such that $\xi_i\oo
\tau_0$ appears in the decomposition of $\psi_*(M)$ are $\tau_0\xi_i$
and $\tau_{i}$ and that for those monomials $\xi_i\oo \tau_0$ appears
with coefficient $1$. The later follows immediately from Lemma
\ref{ongen}. To prove the former note that the question of whether or
not $\xi_i\oo \tau_0$ appears in the expression for $\psi_*(M)$ depend
only on the class of $\psi_*(M)$ in $A_{*,*}\oo_{r,l}A_{*,*}/J$ where
$J$ is generated by $\tau_k$ for $k>0$. In this quotient
$\psi_*(\tau_k)=\xi_k\oo\tau_0+\tau_k\oo 1$ and
$\tau_0^2=\tau\xi_1+\rho\tau_0\xi_1$. This shows that the only way to
get $\xi_i\oo \tau_0$ is to consider $\psi_*(\tau_i)$ or
$\psi_*(\tau_0\xi_i)$.
\end{proof}
The following example shows that not all of the standard topological
formulas for $Q_i$'s hold in the motivic context. 
\begin{example}
\llabel{incor}\rm
In topology, one can define operations $Q_i$ in terms of the Steenrod
squares inductively by the formula $Q_0=Sq^1$,
$Q_{i+1}=[Q_i,Sq^{2^{i+1}}]$. Let us show that in the motivic Steenrod
algebra $Q_2\ne [Q_1,Sq^4]$ if $\rho\ne 0$ i.e. if $k$ does not
contain the square root of $-1$. Using (\ref{copr3}) and Lemma
\ref{ongen} we can compute $Sq^4Q_1$ in terms of the basis dual to
$\tau(E)\xi(R)$. We get:
$$Sq^4Q_1-Q_1Sq^4=Q_2+\rho Q_0Q_1Sq^2$$
\end{example}

\comment{

\subsection{Power operations and topological realization}
In this section we assume that $k={\bf C}$ is the field of complex
numbers. For a pointed simplicial sheaf $F_{\bullet}$, the pointed
simplicial sheaf $G_{\bullet}F_{\bullet}$ has as its terms coproducts
of representable sheaves. Since the category of topological spaces has
coproducts, the simplicial set $G_{\bullet}F_{\bullet}({\bf C})$ has a
natural structure of a simplicial topological spaces. Let
$|G_{\bullet}F_{\bullet}({\bf C})|$ be its geometric realization. This
construction gives us a functor from pointed simplicial sheaves on
$(Sm/\cc)$ to pointed topological spaces $Top_{\bullet}$ which we
denote $t_{\cc}$ and call the topological realization functor.

If $F$ is a coproduct of representable sheaves, then the map
$G_{\bullet}F\sr F$ is a homotopy equivalence. This implies that for
any simplicial pointed sheaf $F_{\bullet}$ whose terms are coproducts
of representables one has a natural homotopy equivalence
$t_{\cc}(F_{\bullet})\sr F_{\bullet}(\cc)$, where the right hand side
is understood as a simplicial topological space. 
\begin{proposition}
\llabel{a1totop}
The functor $t_{\bf C}$ takes $\af$-weak equivalences to homotopy
equivalences. 
\end{proposition}  
\begin{proof}
As shows in \cite{HH2} or in \cite{delnotes}, an $\af$-weak
equivalence between two objects which are term-wise coproducts of
representables belongs to the class $cl_{\bdl}(\{X\times\af\sr X\},
\{K_Q\sr X\})$ which $X$ runs through all smooth schemes and $Q$ runs
through all upper distinguished squares over $X$. Since geometric
realization of simplicial topological spaces commutes with diagonals,
coproducts and filtering colimits it remains to note that morphisms of
the form $(X\times\af)({\bf C})\sr X(\cc)$ and $|K_Q(\cc)|\sr X(\cc)$
are homotopy equivalences of spaces.
\end{proof} 
In view of Proposition \ref{a1totop} the functor $t_{\cc}$ defines a
functor $H_{\bullet}(\cc)\sr H_{\bullet}^{top}$ which we also denote
by $\cc$. 
\begin{proposition}
\llabel{coefcase} For any abelian group $A$, one has isomorphisms in
$H^{top}_{\bullet}$ of the form
$$t_{\bf C}K_{n,A}\sr K(A,2n)$$
These isomorphisms are natural with respect to $A$. 
\end{proposition}
}

\subsection{Operations and characteristic classes}
The goal of this section is to prove Theorem \ref{mainop}. For a
smooth scheme $X$, let $K_0(X)$ be the Grothendieck group of vector
bundles on $X$. All schemes in this section are assumed to be
quasi-projective.
\begin{theorem}
\llabel{char} For any symmetric polynomial $s=s(t_1,\dots,t_n,\dots)$
there exists a unique natural transformation of contravariant functors
from smooth quasi-projective varieties to pointed sets of the form:
$$c_s:K_0(-)\sr \oplus_{n\ge 0} H^{2n,n}(-,\zz)$$
such that for a collection of line bundles $L_1,\dots, L_n$ on $X$ one
has
$$c_s(\oplus_{i=1}^n L_i)=s(e(L_1),\dots,e(L_n))$$
\end{theorem}
\begin{proof}
It follows in the standard way from Theorem \ref{projbund}.
\end{proof}
Let $\phi\in A^{p,q}$ be a cohomological operation. For any $X$ and a
vector bundle $V$ on $X$ the value of $\phi$ on the Thom class $t_V$
is, by Proposition \ref{thom}, of the form $c_{\phi}(V) t_V$
where $c_{\phi}(V)$ is a well defined class in $H^{p,q}(X,\zz/l)$. 
\begin{theorem}
\llabel{mainop}
Let $E=(\epsilon_0,\dots,\epsilon_d)$, $R=(r_1,\dots,r_n)$ be as
in Section \ref{sec13}. Then for a vector bundle $V$ one has:
\begin{enumerate}
\item $c_{\rho(E,R)}(V)=0$ if $E\ne 0$ 
\item $c_{{\cal P}(R)}(V)=s_R(V)$ where $s_R$ is the
reduction modulo $l$ of the characteristic class corresponding by
Theorem \ref{char} to the symmetric polynomial
$$\sum_{f}\prod_{j\in\{0,\dots,n\}}(\prod_{i\in f^{-1}(j)}t_i)^{l^j-1}$$
\end{enumerate}
where $f$ runs through all functions $\{1,\dots,m,\dots\}\sr
\{0,\dots,n\}$ such that for any $i$, $n\ge i>0$ one has
$|f^{-1}(i)|=r_i$.
\end{theorem}
\begin{cor}
\llabel{qicase}
Let $q_n={\cal P}(0,\dots,0,1,0,\dots)$ be the operation dual to
$\xi_n$. Then one has
$$q_n(t_V)=s_{l^n-1}(V) t_V$$
where $s_j$ is the characteristic class corresponding, by Theorem
\ref{char}, to the symmetric function $\sum t^j_i$.
\end{cor}
\begin{cor}
\llabel{powercase}
One has
$$P^n(t_V)=c_{n,l-1}(V) t_V$$
where $c_{n,j}$ is the characteristic class corresponding, by Theorem
\ref{char}, to the symmetric function $\sum t^j_{i_1}\dots t^j_{i_n}$.
\end{cor}
The proof of this theorem occupies the rest of this section.
\begin{lemma}
\llabel{l1.25}
Let $X$ be a smooth scheme and $w$ an element in $H^{2,1}(X,\zz/l)$
which is the reduction modulo $l$ of an integral class. Then there
exists a map $f:X_+\sr (B\mu_l)_+$ in $H_{\BB}$ such that $w=f^*(v)$.
\end{lemma}
\begin{proof}
Since $X$ is quasi-projective the Jouanolou  trick (see \cite{}) shows
that there exist an affine scheme $X'$ and an $\af$-weak equivalence
$X'\sr X$. Therefore, we may assume that $X$ is affine. By \cite{}[]
and Lemma \ref{line}, any element of $H^{2,1}(X,\zz)$ is of the form
$e(L)$ for a line bundle $L$. Since $X$ is affine there is a map $g:X\sr
{\bf P}^{N}$ for some $N$ such that $L=g^*({\cal O}(1))$. On the other
hand the reduction of $e({\cal O}(1))$ modulo $l$ is $p^*(v)$ where
$p:{\bf P}^N\sr B\mu_l$ is the standard morphism. This proves the
lemma.
\end{proof}
\begin{lemma}
\label{l2.25}
Let $X$ be a smooth scheme and $w$ an element in $H^{2,1}(X,\zz/l)$
which is the reduction modulo $l$ of an integral class. Let further
$\phi$ be an operation of the form $\rho(E,R)$. Then one has
$$
\phi(w)=\left\{
\begin{array}{ll}
w^{l^n}&\mbox{\rm if $\phi=q_n$}\\
0&\mbox{\rm otherwise}
\end{array}
\right.
$$
\end{lemma}
\begin{proof}
By Lemma \ref{l1.25} it is sufficient to prove our statement for
$X=B\mu_l$ and $w=v$. In this case our result follows from
(\ref{eqtoth2.2}).
\end{proof}
\begin{lemma}
\llabel{sigpower}
Let $L$ be a line bundle and $\sigma$ the class of ${\cal O}(-1)$ in
$H^{*,*}({\bf P}(L\oplus {\cal O}),\zz)$. Then one has
$\sigma^2=-e(L)\sigma$.
\end{lemma}
\begin{proof}
Using standard argument we can reduce the problem to the case $X={\bf
P}^N$ and $L={\cal O}(1)$. The restriction of $\sigma$ to ${\bf
P}({\cal O})$ is zero. Hence, $\sigma^2=x\sigma$ for some $x$. The
restriction of $\sigma$ to ${\bf P}(L)$ is $-e(L)$. Hence
$-e(L)x=e(L)^2$. Since $e({\cal O}(1))$ is not a zero divisor, we
conclude that $x=-e(L)$.
\end{proof}
\begin{lemma}
\llabel{ofthomnew} Let $L$ be a line bundle and $\phi$ an operation of
the form $\rho(E,R)$. Then one has
$$\phi(t_L)=\left\{
\begin{array}{ll}
e(L)^{l^n-1} t_L&\mbox{\rm for $\phi=q_n$}\\
0&\mbox{\rm otherwise}
\end{array}
\right.$$
\end{lemma}
\begin{proof}
Consider the standard projection $p:{\bf P}({\cal O}\oplus L)\sr
Th(L)$. As shown in Section \ref{sec4}, it defines a monomorphism on
motivic cohomology. Together with Lemma \ref{l2.25} this immediately
implies that $\phi(t_L)=0$ if $\phi\ne q_n$ for some $n$. As shown in
the proof of Lemma \ref{line} we have $p^*(t_L)=\sigma+e(L)$. Hence,
by Lemma \ref{l2.25} and Lemma \ref{sigpower} we have
$$p^*q_n(t_L)=q_np^*(t_L)=(\sigma+e(L))^{l^n}=e(L)^{l^n-1}(\sigma+e(L))$$
Since $p^*$ is a monomorphism we conclude that
$q_n(t_L)=e(L)^{l^n-1} t_L$.
\end{proof}
\begin{remark}
\llabel{adfo}\rm Lemma \ref{ofthomnew} has the following analog for
the basis of admissible monomials instead of Milnor's basis
$\rho(E,R)$. Recall, that $M_k$ denotes the monomial $P^{l^{k-1}}\dots
P^lP^1$. If $L$ is a line bundle and $\phi$ an admissible
monomial then one has
$$\phi(t_L)=\left\{
\begin{array}{ll}
e(L)^{l^k-1} t_L&\mbox{\rm for $\phi=M_k$}\\
0&\mbox{\rm otherwise}
\end{array}
\right.$$
\end{remark}
Let $I$ be the two-sided ideal of $A^{*,*}$ generated by
$Q_0=\beta$. Proposition \ref{com1} implies that it coincides with the
two-sided ideal generated by $Q(E)$ for $E\ne 0$. Since
$\psi^*(\beta)=\beta\oo 1+ 1\oo\beta$ for any $\phi\in I$ we have
$\psi^*(\phi)\in A^{*,*}\oo I + I\oo A^{*,*}$. In particular, if $w,w'$
are motivic cohomology classes such that $\phi(w)=0$, $\phi(w')=0$ for
any $\phi\in I$ then $\phi(w w')=0$ for any $\phi\in I$. Together
with the splitting principle and Lemma \ref{ofthomnew} this implies the
following result.
\begin{lemma}
\llabel{l3.25}
For any $E\ne 0$ and any vector bundle $V$ one has
$c_{\rho(E,R)}(V)=0$.
\end{lemma}
Let $R=(r_1,\dots,r_n)$ be a sequence of non-negative integers. To
prove the second statement of Theorem \ref{mainop} we have to compute
${\cal P}(R)(t_{L_1}\wedge t_{L_m})$ for a collection of line bundles
$L_1,\dots, L_m$. Let $\psi^*_m$ be the m-fold iteration of the
comultiplication map for $A^{*,*}$ and
$$\psi^*_m({\cal P}(R))=\sum a_{(E_1,R_1,\dots, E_m,R_m)}
\rho(E_1,R_1)\oo\dots\oo \rho(E_m,R_m)$$
By Lemma \ref{ofthomnew} we have
$${\cal P}(R)(t_{L_1}\wedge \dots\wedge t_{L_m})=(\sum
a_{(R_1,\dots,R_m)}\prod_{i=1}^m c_{{\cal P}(R_i)}(L_i))\wedge
t_{L_1}\wedge\dots\wedge t_{L_m}$$
where the only non-trivial terms are those for which $R_i$ is of the
form 
$$R_i=(0,\dots,0,1,0,\dots)$$
On the other hand we have
$$a_{(R_1,\dots,R_m)}=\langle \psi^*_m({\cal P}(R)),
\xi(R_1)\oo\dots\oo\xi(R_m)\rangle =$$
$$=\langle  {\cal P}(R),
\xi(R_1)\dots\xi(R_n)=\left\{
\begin{array}{ll}
1&\mbox{\rm if $\sum R_i=R$}\\
0&\mbox{\rm otherwise}
\end{array}
\right.
$$
A sequence of $R_i$'s of the form $(0,\dots,0,1,0,\dots)$ can be
thought of as a function $f:i\sr f(i)$ such that $R_i=q_{f(i)}$ where
$q_0$ is assumed to be $1$. The condition $\sum R_i=R$ means that we
consider the functions which take the value $i>0$ exactly $r_i$
times. Together with the fact that $c_{q_n}(L)=e(L)^{l^n-1}$ this
proves the last statement of Theorem \ref{mainop}.


\begin{thebibliography}{10}

\bibitem{SBloch1}
S.~Bloch.
\newblock The moving lemma for higher {C}how groups.
\newblock {\em J. Algebr. Geom.}, 3(3):537--568, Feb. 1994.

\bibitem{Boardman}
J.~M. Boardman.
\newblock The eightfold way to {BP}-operations.
\newblock In {\em Current trends in algebraic topology}, pages 187--226.
  AMS/CMS, Providence, 1982.

\bibitem{delnotes}
Pierre Deligne.
\newblock Voevodsky's lectures on motivic cohomology.
\newblock {\em www.math.ias.edu/{$_{\textstyle
  \tilde{}}\,$}vladimir/seminar.html}, 2000/2001.

\bibitem{Milnor}
J.~Milnor.
\newblock Algebraic {K}-theory and quadratic forms.
\newblock {\em Inv. Math.}, 9:318--344, 1970.

\bibitem{Milnor3}
John Milnor.
\newblock The {S}teenrod algebra and its dual.
\newblock {\em Annals of Math.}, 67(1):150--171, 1958.

\bibitem{MoVo}
Fabien Morel and Vladimir Voevodsky.
\newblock {${\bf A}^1$}-homotopy theory of schemes.
\newblock {\em Publ. Math. IHES}, (90):45--143, 1999.

\bibitem{SE}
N.~E. Steenrod and D.~B. Epstein.
\newblock {\em Cohomology operations}.
\newblock Princeton Univ. Press, Princeton, 1962.

\bibitem{SusVoe3}
Andrei Suslin and Vladimir Voevodsky.
\newblock {B}loch-{K}ato conjecture and motivic cohomology with finite
  coefficients.
\newblock In {\em The arithmetic and geometry of algebraic cycles}, pages
  117--189. Kluwer, 2000.

\bibitem{MC}
Vladimir Voevodsky.
\newblock The {M}ilnor {C}onjecture.
\newblock {\em www.math.uiuc.edu/K-theory/170}, 1996.

\bibitem{comparison}
Vladimir Voevodsky.
\newblock Motivic cohomology are isomorphic to higher {C}how groups.
\newblock {\em www.math.uiuc.edu/K-theory/378}, 1999.

\bibitem{MCnew}
Vladimir Voevodsky.
\newblock On 2-torsion in motivic cohomology.
\newblock {\em Preprint}, 2001.

\end{thebibliography}
\end{document}